\documentclass[11pt,a4paper]{article}
\usepackage{amsmath}
\usepackage{amssymb}
\usepackage{color}

\usepackage{graphicx}
\usepackage{amsfonts,amsmath, amssymb}
\usepackage{hyperref}

\numberwithin{equation}{section}

\usepackage{graphics}
\usepackage{epsfig}
\newtheorem{claim}{\bf \t}[part]

\newtheorem{theorem}{Theorem}[section]
\newtheorem{corollary}[theorem]{Corollary}
\newtheorem{lemma}[theorem]{Lemma}
\newtheorem{proposition}[theorem]{Proposition}
\newtheorem{remark}[theorem]{Remark}

\def\v{\varepsilon}
\def\e{\epsilon}

\def\t{\theta}
\def\k{\kappa}

\def\m{\mu}
\def\a{\alpha}
\def\b{\beta}
\def\g{\gamma}
\def\d{\delta}
\def\l{\lambda}

\def\r{\rho}
\def\s{\sigma}

\def\i{\infty}
\def\f{\frac}
\begin{document}

	\title{Diffusive Wave in the Low Mach  Limit for  Compressible Navier-Stokes Equations}
	
	\author{Feimin Huang$^{a}$, Tian-Yi Wang$^{b,c}$,
		Yong Wang$^{a}$\footnote{
			Email addresses: fhuang@amt.ac.cn(Feimin Huang), tian-yi.wang@gssi.infn.it(Tian-Yi Wang),  yongwang@amss.ac.cn(Yong Wang)}
		\\
		\ \\
		{\small \it $^a$Institute of Applied Mathematics, AMSS, CAS, Beijing 100190, China}\\[1mm]
		{\small \it $^b$	Department of Mathematics, School of Science, Wuhan University of Technology, Wuhan,  China}\\[1mm]
		{\small \it $^c$Gran Sasso Science Institute, viale Francesco Crispi, 7, 67100 L'Aquila, Italy}
		}
	
	\date{ }
	\maketitle
	
	\begin{abstract}
		
		The low Mach limit for 1D non-isentropic compressible Navier-Stokes flow, whose density and temperature have different asymptotic states at infinity, is rigorously justified. The problems are considered on  both well-prepared and ill-prepared data. For the well-prepared  data, the solutions of compressible Navier-Stokes equations are shown to converge to a nonlinear diffusion wave solution globally in time as Mach number goes to zero when the difference between the states at $\pm\infty$ is suitably small. In particular, the velocity of diffusion wave is only driven by the variation of temperature.  It is  further shown that the solution of  compressible Navier-Stokes system also has the same property when Mach number is small, which has never been observed before. The convergence rates on both Mach number and time are also obtained for the well-prepared data. For the ill-prepared data, the limit relies on the uniform estimates including  weighted time derivatives and an extended convergence lemma.  And the difference between the states at $\pm\infty$  can be arbitrary large  in this case.
		
		\
		
		Keywords: Compressible Navier-Stokes equations, low Mach limit, nonlinear diffusion wave, well-prepared data, ill-prepared data
		
		\
		
		AMS: 35Q35, 35B65, 76N10 
	\end{abstract}
	
	\tableofcontents
	
	\section{Introduction}
	
	The non-isentropic Navier-Stokes system in $\mathbb{R}^n$ is as follows
	\begin{align}\label{1.1}
		\begin{cases}
			\partial_t\r+\mbox{div}(\r u)=0,\\
			\partial_t(\r u)+\mbox{div}(\r u\otimes u)+\nabla P=\mbox{div}(2\mu D(u))+\nabla(\l \mbox{div}u),\\
			\partial_t(\r(e+\f12|u|^2))+\mbox{div}(\r u(e+\f12|u|^2)+P u)=\mbox{div}(\k \nabla \mathcal{T})+\mbox{div}(2\mu D(u) u+\l \mbox{div}u u), 
		\end{cases}
	\end{align}
	for $t>0, x\in \mathbb{R}^n$. Here the unknown functions $\r$, $u$, and $\mathcal{T}$ represent the density, velocity, and temperature, respectively. The pressure function and internal function are defined  by
	\begin{equation}\label{1.6}
		P=R\r\mathcal{T},~~e=c_v\mathcal{T},
	\end{equation}
	where the parameters  $R>0$ and $c_v>0$ are the gas constant and the heat capacity at constant volume, respectively.   $D(u)$ is the deformation tensor given by
	$$
	D(u)=\f12(\nabla u+(\nabla u)^{t}),
	$$
where $(\nabla u)^t$ denote the transpose of matrix $\nabla u$. $\m$ and $\l$ are the Lam\'{e} viscosity coefficients which satisfy
$$\mu>0,~~2\mu+n\l>0,$$
and $\k>0$ is the heat conductivity coefficient. For simplicity, we assume that  $\mu,~\l,~\k$ are constants.

It is noted $c_v=(\g-1)/R$ where $\g>1$ is the adiabatic exponent. In this paper, since we consider the low Mach limit for smooth solutions of non-isentropic compressible Navier-Stokes equation,  we normalize  $R=1$ and $c_v=1$ for simplicity of presentation, see also \cite{JJLX}. And we point out that our results and energy estimates in this paper still hold for any given $c_v>0$(equivalently for any given $\g>1$), the only difference is that some constants of this paper may depend on $c_v$.
	
Let $\v$ be the compressibility parameter which is a nondimensional quantity. As in \cite{Schochet2}, we set
\begin{equation}\label{1.2}
t\rightarrow \v t,~~x\rightarrow x,~~u\rightarrow \v u,~~ \mu\rightarrow \v\mu,~~\l\rightarrow \v\l,~~\k\rightarrow \v\k.
\end{equation}
Based on the above changes of variables, the compressible Navier-Stokes system \eqref{1.1}, written after the nondimensionalization, becomes
\begin{eqnarray}\label{1.3}
\begin{cases}
\partial_t\r^\v+\mbox{div}(\r^\v u^\v)=0,\\
\partial_t(\r^\v u^\v)+\mbox{div}(\r^\v u^\v\otimes u^\v)+\f{\nabla P^\v}{\v^2}=\mbox{div}(2\mu D(u^\v))+\nabla(\l \mbox{div}u^\v),\\
\partial_t(\r^\v\mathcal{T}^\v)+\mbox{div}(\r^\v u^\v\mathcal{T}^\v)+P^\v\mbox{div}u^\v=\mbox{div}(\k \nabla\mathcal{T}^\v)+\v^2[2\mu |D(u^\v)|^2 +\l |\mbox{div}u^\v|^2],
\end{cases}
\end{eqnarray}
in which the typical mean fluid velocity has been chosen as the ratio of time units to space units. Accordingly, the parameter $\varepsilon$ essentially presents the maximum Mach number of the fluid. The limit of solutions of \eqref{1.3} as   $\v$ goes to $0$ is usually called as low Mach limit \cite{KM1, KM2}.  The purpose of the low Mach number approximation is to justify that the compression, due to pressure variations, can be neglected. This is a common assumption  when discussing the fluid dynamics of highly subsonic flows. Similar to \cite{Alazard2,Schochet2}, we assume that the pressure is a small perturbation of a given constant state $\underline{P}>0$, \textit{i.e.}
\begin{equation}\label{1.4}
P^\varepsilon=\underline{P}+O(\v), 
\end{equation}
which satisfies the following equation
\begin{equation}\label{1.3+}
\partial_t(P^\v)+\mbox{div}( u^\v P^\v)+P^\v\mbox{div}u^\v=\mbox{div}(\k \nabla\mathcal{T}^\v)+\v^2[2\mu |D(u^\v)|^2 +\l |\mbox{div}u^\v|^2].
\end{equation}
Formally, as the Mach number tends to zero, the equations \eqref{1.3+}, $(\ref{1.3})_2$ and $(\ref{1.3})_3$ become
\begin{equation}\label{1.5}
\begin{cases}
\mbox{div}(2\underline{P}u-\k\nabla\mathcal{T})=0,\\
\r (u_t+u\cdot\nabla u)+\nabla \pi=\mbox{div}(2\mu D(u))+\nabla(\l \mbox{div}u),\\
\r (\mathcal{T}_t+u\cdot\nabla\mathcal{T})+\underline{P}\mbox{div}u=\mbox{div}(\k\nabla\mathcal{T}),
\end{cases}
\end{equation}
and the density satisfies $\r=\underline{P}/\mathcal{T}$ and 
\begin{equation}
\r_t+\mbox{div}(\r u)=0.
\end{equation}
Without loss of generality, we  normalize $\underline{P}$ to be $1$. Let the approximate initial data $(p^{\v}_{in}, u^{\v}_{in}, \mathcal{T}^{\v}_{in})$ converge to $({p}_{in}, {u}_{in}, {\mathcal{T}}_{in})$ as $\v\rightarrow0$ in some sense.  Then, the approximate initial data is regarded as well-prepared data if  ${u}_{in}$ and ${\mathcal{T}}_{in}$ satisfies  $(\ref{1.5})_1$. Otherwise, it is  called as  ill-prepared data.
	
The low Mach limit is an important and interesting  problem in  the fluid dynamics. There have been a lot of literatures on the low Mach Limit. The first result is due to Klainerman and Majda \cite{KM1, KM2}, in which they proved the incompressible limit of the isentropic  Euler equations to the incompressible Euler equations for  local smooth solutions with  well-prepared  data.  By using the fast decay of acoustic waves, Ukai \cite{Ukai} verified the low Mach limit for the general  data, see also \cite{Asano}. The half plane  and exterior domain cases were considered in  \cite{Iguchi, Isozaki1, Isozaki2}. In \cite{Schochet1}, Schochet showed the limit of the full compressible Euler equations to the incompressible inhomogeneous Euler equations in bounded domain for local smooth solutions with well-prepared initial data. He \cite{Schochet}  further studied the fast singular limit for general hyperbolic partial differential equations.  The major breakthrough on the ill-prepared data  is due to M\'{e}tivier and Schochet \cite{MS}, in which they proved the low Mach limit of full Euler equations in the whole space by a significant convergence lemma on acoustic waves. Later, Alazard \cite{Alazard1} extended \cite{MS} to the exterior domain. The one dimension spatial periodic case on the full compressible Euler equations was considered in \cite{MS1}. For other interesting works,  see \cite{BDG} and the references therein.

For the isentropic Navier-Stokes equations, the low Mach limit of the global weak solutions with general initial data have been well studied under  various boundary conditions, see \cite{ DeG, DGLM, Lions-Masmoudi, Masmoudi}. For the non-isentropic  Navier-Stokes equations, Alazard \cite{Alazard2} justified  the low Mach limit in the whole space for the ill-prepared data, by employing a uniform estimate and the convergence lemma of \cite{MS}.  For the bounded domain, the low Mach limit was justified by Jiang-Ou \cite{Jiang-Ou} and Dou-Jiang-Ou \cite{Dou-Jiang-Ou}.  A dispersive Navier-Stokes system was also studied in \cite{Levermore}.  For other interesting works, see  \cite{DBDL,Danchin,Feireisl-N,Kim,Masmoudi1,Schochet2} for Naiver-Stokes equations, \cite{Fan-Gao-Guo,Hu-Wang,JJL1,JJL2,JJL3, JJLX} for  MHD equations and the references therein.

Note that in the whole space, all results above require that both density and temperature have the constant background, \textit{i.e.}
\begin{equation}\label{1.7}
\r^\v(x,t)\rightarrow \underline{\r},~\mbox{and}~\mathcal{T}^\v(x,t)\rightarrow \underline{\mathcal{T}},~~\mbox{as}~~|x|\rightarrow \infty,
\end{equation}
where $\underline{\r}>0$ and $\underline{\mathcal{T}}>0$ are given constants. In this paper, we will study the low Mach limit when the background is not constant state in the one dimensional case, that is
\begin{equation}\label{1.8}
(\r^\v,\mathcal{T}^\v)(x,t)\rightarrow (\r_{\pm},\mathcal{T}_{\pm}), ~~\mbox{as}~x\rightarrow \pm\infty, ~~ \mbox{with}~~\r_-\mathcal{T}_-=\r_+\mathcal{T}_+,
\end{equation}
where $\mathcal{T}_-$ may not be equal to $\mathcal{T}_+$, and want to know what happens in the limiting process. In fact, we find the solutions of compressible Navier-Stokes equations converge to a nonlinear diffusion wave solution globally in time as Mach number goes to zero. In particular, the velocity of diffusion wave is only driven by the variation of temperature.  This phenomenon looks like thermal creep flow \cite{Huang,Huang-Wang-Wang-Yang-2,Slemrod2,Slemrod}. Moreover,  when Mach number is small, the compressible Navier-Stokes system also has the same property, which has never been observed before. 

Now we begin to formulate the main results.  The system \eqref{1.3} and  the limiting system \eqref{1.5} in one dimension become, respectively, 
\begin{align}\label{1.3-1}
\begin{cases}
\partial_t\r^\v+(\r^\v u^\v)_x=0,\\
\r^\v (u^\v_t+ u^\v u^\v_x)+\f{ P^\v_x}{\v^2}=\tilde{\m}u^{\v}_{xx},\\
\r^\v(\mathcal{T}^\v_t+ u^\v\mathcal{T}^\v_x)+P^\v u^\v_{x}=\k\mathcal{T}^\v_{xx}+\tilde{\m}\v^2|u^\v_x|^2,
\end{cases}
\end{align}
and
\begin{equation}\label{1.5-1}
\begin{cases}
(2u-\k\mathcal{T}_x)_x=0,~~\r=\mathcal{T}^{-1},\\
\r (u_t+u u_x)+\pi_x=\tilde{\m}u_{xx},\\
\r (\mathcal{T}_t+u\mathcal{T}_x)+u_x=\k\mathcal{T}_{xx},
\end{cases}
\end{equation}
where $(x,t)\in \mathbb{R}\times\mathbb{R}_+$ and  $\tilde{\m}=2\mu+\l>0$. 
We shall construct a special solution of \eqref{1.5-1}. Indeed, from $\eqref{1.5-1}_1$, we choose
\begin{align}\label{1.11}
u=\f{\k}{2}\mathcal{T}_x=-\f{\k}{2}\f{\r_x}{\r^2}.
\end{align}
Substituting \eqref{1.11} into \eqref{1.5-1}, one obtains the following nonlinear diffusion equation
\begin{align}\label{1.12}
\r_t=\Big(\f{\k}{2}\f{\r_x}{\r}\Big)_x.
\end{align}
We consider the boundary condition at the far field, i.e.,  $\lim_{x\rightarrow\pm\infty}\r(x,t)=\r_{\pm}$. 
From \cite{Atkinson-Peletier} and \cite{Duyn-Peletier}, it is known that
the nonlinear diffusion equation \eqref{1.12} admits a unique
self-similar solution $\Xi(\eta),~\eta=\f{x}{\sqrt{1+t}}$
satisfying $
\Xi(\pm\i,t)=\r_\pm. $ Let $\delta=|\r_+-\r_-|$, 
then $\Xi(t,x)$ satisfies
\begin{equation*}
\Xi_x(t,x)=\f{O(1)\delta}{\sqrt{(1+t)}}e^{-\f{x^2}{4d(\r_{\pm})(1+t)}}, ~~\mbox{as}~~x\rightarrow\pm\infty, ~~\mbox{where} ~~ d(\r)=\f{\k}{2\r}.
\end{equation*}
We define 
\begin{align}\label{1.15}
(\bar{\rho},\bar{u},\bar{\mathcal{T}})\doteq \left(\Xi, ~-\f{\k}{2}\f{\Xi_x}{\Xi^2},~\Xi^{-1}\right),
\end{align}
which is a special solution of \eqref{1.5-1}, that is
\begin{equation}\label{1.5-2}
\begin{cases}
(2\bar u-\k\bar{\mathcal{T}}_x)_x=0,~\bar{\r}=\bar{\mathcal{T}}^{-1},\\
\bar\r (\bar u_t+\bar u \bar u_x)+\bar{\pi}_x=\tilde{\m}\bar{u}_{xx},\\
\bar\r (\bar{\mathcal{T}}_t+\bar u\bar{\mathcal{T}}_x)+\bar{u}_x=\k\bar{\mathcal{T}}_{xx},
\end{cases}
\end{equation}
with $\bar{\pi}=\tilde{\m}\bar{u}_{x}-\bar{\rho}\bar{u}^2+\f{\k}{2}\f{\bar{\r}_t}{\bar{\r}}$. 
In particular, $\bar{u}=\frac{\kappa}{2}\bar{\mathcal{T}}_x$ indicates that the flow 
is driven by the variation of temperature. In other word, the flow moves 
along the direction from the low temperature to the high one, see Figures 1 and 2 below, and thus behaves like thermal creep flow.
\begin{figure}[htbp]
	\begin{center}
		\includegraphics[width=12cm]{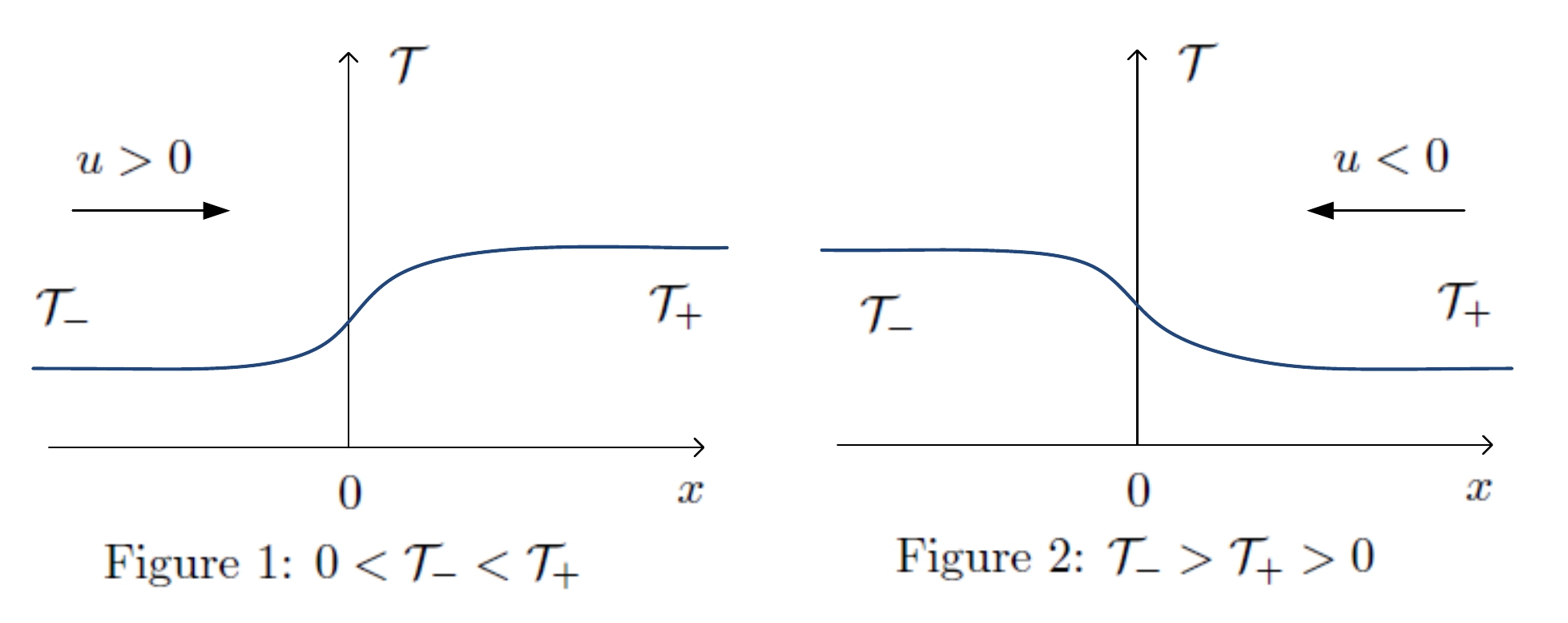}
	\end{center}
\end{figure}
Indeed, Slemrod \cite{Slemrod}  mentioned:  \textquotedblleft{\textit{This leaves open the issue of finding the correct extended fluid dynamics for moderate dense gas or polyatomic gas}\textquotedblright, \textit{i.e.}, how to construct a thermal creep flow for moderate dense gas or polyatomic gas is still an open problem. In this paper,  for well-prepared data, we will show that the solution of compressible Navier-Stokes equation behaviors as thermal creep flow in some sense when the Mach number is small. It should be noted that the usual Poiseuille flow moves  from the high pressure part to the low one due to the difference of pressure.  We will state our main results in the following two subsections.

\subsection{Main Results for Well-prepared Data}
In  this subsection, we consider the low Mach number limit around the diffusive wave with well-prepared initial data. It is more convenient to use the {\it Lagrangian} coordinates for the global behavior of solutions. That is, take the coordinate transformation
\begin{equation}\nonumber
	(x,t)\rightarrow\Big(\int_{(0,0)}^{(x,t)}\r^\v(z,s)dz-(\r^\v u^\v)(z,s)ds, s\Big),
\end{equation}
which is  still denoted  as $(x,t)$ without confusion. It is noted that the above transformation is independent of the path of integration due to the continuity equation $\eqref{1.3-1}_1$.   We also point out that  the Eulerian and  Lagrangian formulations are equivalent since we consider the smooth solution of compressible Navier-Stokes equations, see also \cite{Batchelor,Serrin}.  Let $v=\f{1}{\r}$ be the specific volume.  Then, in  Lagrangian coordinates,   the system \eqref{1.3-1}  can be written equivalently as
\begin{eqnarray} \label{N-1.1}
	\begin{cases}
		v^\v_{t}-u^\v_{x}=0, \\[2mm]
		u^\v_{t}+\f1{\v^2} P^\v_{x}=\tilde{\m}(\frac{u^\v_x}{ v^\v})_{x},\\[2mm]
		(\mathcal{T}^\v+\f12|\v u^\v|^2)_t+(P^\v u^\v)_x=\k(\frac{1}{v^\v}\mathcal{T}^\v_x)_x+\v^2(\frac{\tilde\mu}{v^\v}u^\v u^\v_{x})_x,
	\end{cases}
\end{eqnarray}
where the pressure $P=\mathcal{T}/v$. 
Similarly, the limiting system \eqref{1.5-2} in the Lagrangian coordinate becomes
\begin{equation}\label{N-1.5-1}
	\begin{cases}
		(2u-\f{\k}{v}\mathcal{T}_x)_x=0,~v=\mathcal{T},\\
		u_t+\pi_x=\tilde{\m}(\f{1}{v}u_{x})_x,\\
		\mathcal{T}_t+ u_x=\k(\f1{v}\mathcal{T}_x)_x.
	\end{cases}
\end{equation}
As in \eqref{1.15}, we can construct a special diffusive wave solution $(\bar{v},\bar{u},\bar{\mathcal{T}},\bar{\pi})$ of \eqref{N-1.5-1} by choosing 
\begin{equation}\label{N-3.6}
(\bar{v},\bar{u},\bar{\mathcal{T}}) := (\hat{\mathcal{T}}, \f{\k\hat{\mathcal{T}}_x}{2\hat{\mathcal{T}}},\hat{\mathcal{T}}),
\end{equation}
where $\hat{\mathcal{T}}(\eta), \eta=\f{x}{\sqrt{1+t}}$ is the unique self-similar solution of the following diffusion equation 
\begin{equation}\label{N-1.3}
\mathcal{T}_t=\left(\f{\k\mathcal{T}_x}{2\mathcal{T}}\right)_x,~\mbox{with}~\lim_{x\rightarrow\pm\infty}\mathcal{T}=\mathcal{T}_{\pm},
\end{equation}
and 
$\bar{\pi}$ can be solved by $\eqref{N-1.5-1}_2$.
Let $\delta=|\mathcal{T}_+-\mathcal{T}_-|$, then
$\hat{\mathcal{T}}(t,x)$ satisfies
\begin{equation}\label{N-1.4}
	\hat{\mathcal{T}}_x(t,x)=O(1)\delta(1+t)^{-\f12}e^{-\f{x^2}{4d(\mathcal{T}_{\pm})(1+t)}}, ~~d(\mathcal{T})=\f{\k}{2\mathcal{T}},~~\mbox{as}~~ x\rightarrow\pm\infty.
\end{equation}
Indeed, $(\bar{v},\bar{u},\bar{\mathcal{T}})$ is the Lagrangian version of diffusion wave solution  \eqref{1.15}.

Since $(\bar{v},\bar{u},\bar{\mathcal{T}})$ is not the solution of \eqref{N-1.1} and some non-integrated error terms with respect to time $t$ and slow  $\v$-decay terms should appear in the low Mach limiting analysis, we need to introduce a new profile to approximate the system \eqref{N-1.1}, \textit{i.e.},  
\begin{equation}\label{N-3.8}
	(\tilde{v},\tilde{u},\tilde{\mathcal{T}})
	:= (\bar{v}, \bar{u}, \bar{\mathcal{T}}-\f{1}{2}|\v\bar{u}|^2).
\end{equation}
Then a direct calculation gives that
\begin{align} \label{N-1.5}
	\begin{cases}
		\tilde v_{t}-\tilde u_{x}=0, \\
		\v^2\tilde u_{t}+\tilde P_{x}=\tilde\mu\v^2(\frac{\tilde u_{x}}{\tilde v})_{x}+\tilde R_{1x},\\
		\left(\tilde{\mathcal{T}}+\frac{1}{2}|\v\tilde u|^{2}\right)_{t}+
		(\tilde P\tilde u)_{x}=(\k\frac{\tilde{\mathcal{T}}_{x}}{\tilde v})_x+\v^2 (\frac{\tilde\mu}{\tilde v}\tilde u\tilde u_{x})_{x}+\tilde R_{2x},
	\end{cases}
\end{align}
where
\begin{align}
	&\tilde R_1=\v^2\f{\k\hat{\mathcal{T}}_t}{2\hat{\mathcal{T}}}-\v^2\f{\tilde{u}^2}{2\tilde{v}}-\v^2\f{\tilde\mu}{\tilde v}\tilde{u}_x=O(1)\d\v^2(1+t)^{-1}e^{-\f{c_\pm x^2}{1+t}},~~\mbox{as}~~|x|\rightarrow\infty,\label{N-1.6}\\
	&\tilde R_2=\v^2\f{\k}{\hat{\mathcal{T}}}\tilde{u}\tilde{u}_x-\f{\v^2}{2\tilde v}\tilde{u}^3
	-\f{\v^2}{\tilde v}\tilde{\mu}\tilde u\tilde u_x=O(1)\d\v^2(1+t)^{-\f32}e^{-\f{c_\pm x^2}{1+t}},~~\mbox{as}~~|x|\rightarrow\infty \label{N-1.7}.
\end{align}
It is noted that the system \eqref{N-1.5} approximate the original compressible Navier-Stokes equations \eqref{N-1.1} very well, up to error terms $\tilde{R}_{1x}$ and $\tilde{R}_{2x}$, when $\v$ is sufficiently small. Moreover, since we will integrate the the differences between the solutions of original system \eqref{N-1.1}  and the approximate system \eqref{N-1.5} with respect to space variable in \eqref{N-2.4} below, so it is very important that the error terms are in the form of $x$-derivatives. From \eqref{N-3.6} and \eqref{N-3.8}, we also note  that 
	\begin{align}\label{N-3.12}
		\|(\bar{v}-\tilde{v},\bar{u}-\tilde{u},\bar{\mathcal{T}}-\tilde{\mathcal{T}})(t)\|\leq C\v^2 (1+t)^{-1},
	\end{align}
which implies that  $(\tilde{v},\tilde{u},\tilde{\mathcal{T}})$ approximate the diffusive wave solution $(\bar{v},\bar{u},\bar{\mathcal{T}})$ very well  when $\v$ is small.

Now we supplement the system \eqref{N-1.1} with  the initial data
\begin{equation}\label{N-1.7-1}
	(v^\v,u^\v,\mathcal{T}^\v)|_{t=0}=(\tilde{v}, \tilde{u}, \tilde{\mathcal{T}})(x,0),
\end{equation}
then we have  the following global existence and uniform estimates.
\begin{theorem}[Uniform Estimates for Well-Prepared Data]\label{Nthm1.1}
	Let $(\tilde{v}, \tilde{u}, \tilde{\mathcal{T}})(x,t)$ be the diffusive wave defined in  \eqref{N-3.8} with the wave strength $\d=|\mathcal{T}_+-\mathcal{T}_-|$. There exist positive constants $\d_0$ and $\v_0$, such that if $\d\leq \d_0$  and $\v\leq \v_0$, then the Cauchy problem  \eqref{N-1.1}, \eqref{N-1.7-1} has a unique global smooth solution $(v^\v,u^\v,\mathcal{T}^\v)$ satisfying
	\begin{equation}\label{N-1.7-3}
		\begin{cases}
			\|(v^\v-\tilde{v},\v u^\v-\v\tilde{u},\mathcal{T}^\v-\tilde{\mathcal{T}})(t)\|^2_{L^2_x}\leq C \sqrt{\d}\v^3(1+t)^{-1+C_0\sqrt{\d}},\\
			\|(v^\v-\tilde{v},\v u^\v-\v\tilde{u},\mathcal{T}^\v-\tilde{\mathcal{T}})_x(t)\|^2_{L^2_x}\leq C \sqrt{\d}\v^2(1+t)^{-\f32+C_0\sqrt{\d}},\\
		\end{cases}
	\end{equation}
	where $C$ and $C_0$ are positive constants independent of $\v$ and $\d$.
\end{theorem}

Based on the uniform estimates \eqref{N-1.7-3} and Sobolev embedding, we justify the following  low Mach limit. 
\begin{corollary}[Low Mach Limit for Well-Prepared Data]\label{Ncor1.2}
	Under the assumptions of Theorem \ref{Nthm1.1}, it follows from \eqref{N-3.12} and \eqref{N-1.7-3}  that  as $\v\rightarrow0$,
	\begin{equation}\label{N-1.7-2}
		\begin{cases}
			\|(v^\v-\bar{v},\mathcal{T}^\v-\bar{\mathcal{T}})(t)\|_{L^\infty({\mathbb R})}\leq C \d^\f14\v^{\f54}(1+t)^{-\f12}
			\rightarrow 0,\\
			\|(u^\v-\bar{u})(t)\|_{L^\infty({\mathbb R})}\leq C \d^\f14\v^{\f14}(1+t)^{-\f12}
			\rightarrow 0.
		\end{cases}
	\end{equation}

\end{corollary}

\

Note that the velocity  $\bar{u}$ is  driven by the variation of temperature $\tilde{\mathcal{T}}$, \textit{i.e.}, $\bar{u}=\f{\k\hat{\mathcal{T}}_x}{2\hat{\mathcal{T}}}$. Without loss of generality, we assume that $\mathcal{T}_- < \mathcal{T}_+$, then for any given positive constant $\eta_0>0$, there exists $c_{\eta_0}>0$ such that 
\begin{eqnarray}\label{N1.8}
	\hat{\mathcal{T}}'(\eta)>c_{\eta_0}\d>0,~~ \mbox{for}~~ |\eta|\leq \eta_0.
\end{eqnarray}
Moreover, we have
\begin{corollary}[Driven by the Variation of Temperature]\label{Ncor1.3}
	For the solutions obtained in Theorem \ref{Nthm1.1}, there exists a small positive constant $\v_1=\v_1(\eta_0)\leq \v_0$, such that if $\v\leq \v_1$,   it follows from \eqref{N-1.7-2} and \eqref{N1.8} that
	\begin{equation}\label{N1.9}
		\begin{cases}
			0<\f{c_{\eta_0}\d}{C_{1}\sqrt{1+t}}<\f{1}{C_{1}}\hat{\mathcal{T}}_x\leq u^\v(x,t)\leq C_{1}\hat{\mathcal{T}}_x,\\[3mm]
			~~~~~~~~0<\f12 \hat{\mathcal{T}}_x\leq \mathcal{T}_x^\v(x,t)\leq \f32\hat{\mathcal{T}}_x,
		\end{cases}
		\mbox{\rm for}~ |x|\leq \eta_0(1+t)^{\f12},~t\geq0,
	\end{equation}
	which yields immediately that 
	\begin{equation}\label{B1.55-1}
		\f{2}{3C_1}\mathcal{T}^\v_x(x,t)\leq u^\v(x,t)\leq 2C_1\mathcal{T}^\v_x(x,t),~\mbox{\rm for}~ |x|\leq \eta_0(1+t)^{\f12},~t\geq 0,
	\end{equation}
	where $C_{1}$ is a suitably large positive constant  depending only on $\mathcal{T}_\pm$.
\end{corollary}

\begin{remark}
	The estimate \eqref{B1.55-1}  shows that the velocity  $u^\v$ of the compressible Navier-Stokes system \eqref{N-1.1} is also driven by the variation of temperature when $\v$ is small. 	Note that the pressure $P$ is almost 1, this phenomenon behaviors like 
	thermal creep flow \cite{Huang,Slemrod2,Slemrod}. So, our result can be  regarded as an answer to the open question of \cite{Slemrod} in some sense.
\end{remark}
\begin{remark}
	All the above results for the well-prepared data hold globally in time. 
\end{remark}

\subsection{Main Results for Ill-prepared Data}
To understand the role of the thermodynamics, as in \cite{Alazard2},   we rewrite the equations \eqref{1.3} by the pressure fluctuations $p^\v$, and temperature fluctuations $\t^\v$ with 
\begin{equation}\label{2.4-1}
	P^\v(x,t)=e^{\v p^\v(x,t)},~~\mathcal{T}^\v(x,t)=e^{\t^\v(x,t)}.
\end{equation}
It follows from \eqref{1.6} and \eqref{2.4-1} that
\begin{equation}\label{2.7}
	\r^\v(x,t)=e^{\v p^\v(x,t)-\t^\v(x,t)}.
\end{equation}
Under these changes of variables, the asymptotic states at infinity are, as $x\rightarrow\pm\infty$
\begin{equation}
	p^\v\rightarrow 0,~~ u^\v\rightarrow 0, ~~\mbox{and} ~~\theta^\v\rightarrow \theta_\pm,~\mbox{as}~x\rightarrow\pm\infty, ~~\mbox{where}~~ \theta_\pm=\ln \mathcal{T}_{\pm}.
\end{equation}
The  compressible Navier-Stokes system \eqref{1.3-1}  takes the following equivalent form
\begin{align}\label{2.1}
	\begin{cases}
		p^\v_t+u^\v\cdot p^\v_x+\f{1}{\v}(2u^\v-\k e^{-\v p^\v+\t^\v}\t^\v_x)_x
		=\tilde\m\v e^{-\v p^\v}|u^\v_x|^2+\k e^{-\v p^\v+\t^\v}p^\v_x\cdot\t^\v_x ,\\[2mm]
		e^{-\t^\v}[u^\v_t+u^\v\cdot u^\v_x]+\f1\v p^\v_x=\tilde{\mu}e^{-\v p^\v}u^\v_{xx},\\[2mm]
		\t^\v_t+u^\v\t^\v_x+u^\v_x=\k e^{-\v p^\v}(e^{\t^\v}\t^\v_x)_x+\tilde{\m}\v^2e^{-\v p^\v}|u^\v_x|^2,
	\end{cases}
\end{align}
and the limiting system \eqref{1.5-1} becomes
\begin{equation}\label{2.8}
	\begin{cases}
		(2u-\k e^{\t}\t_x)_x=0,\\
		e^{-\t}(u_t+u u_x)+\pi_x=(\tilde{\mu}u_x)_x,\\
		\t_t+u \t_x+u_x=(\k e^{\t}\t_x)_x.
	\end{cases}
\end{equation}


Since $\theta_{-}$ may not be equal to $\theta_{+}$, we need to introduce a background profile $	\tilde{\t}$ for $\theta^\v$. Here we choose 
\begin{equation}
	\tilde{\t}:=- \ln \Xi,
\end{equation}
satisfying $\tilde{\theta}\rightarrow\theta_{\pm}$ as $x\rightarrow\pm\infty$.
Then we define the following solution space: for $s\geq 4$,
\begin{align}\label{2.2}
	\mathcal{N}(t)&:=\|(p^\v,u^\v,\t^\v-\tilde{\t})(t)\|^2_{H^{s,\v}}\nonumber\\
	&=\sum_{|\a|=0}^s\|\partial^\a(p^\v,u^\v)(t)\|^2_{L^2}
	+\sum_{|\a|=0}^{s+1}\|\partial^\a(\v p^\v, \v u^\v)(t)\|^2_{L^2}+\|(\t^\v-\tilde{\t})(t)\|^2_{L^2}+\sum_{|\a|=1}^{s+1}\|\partial^\a \t^\v(t)\|^2_{L^2}\nonumber\\
	&~~~~~~+\sum_{|\a|=0}^s\int_0^t\|\partial^\a(p^\v_x,u^\v_x)(\tau)\|^2_{L^2}d\tau
	+\sum_{|\a|=0}^{s+1}\int_0^t\|\partial^\a(\v u^\v_x,\t^\v_x)(\tau)\|^2_{L^2}d\tau,
\end{align}
where $\partial^\a:=(\v\partial_t)^{\a_0}\partial_x^{\a_1}$ with multi-index $\a=(\a_0,\a_1)$.

We  supplement the system \eqref{2.1} with  the following  initial data
\begin{align}\label{2.1-1}
	(p^\v,u^\v,\t^\v)|_{t=0}=(p^\v_{in},u^\v_{in},\t^\v_{in}),
\end{align}
with
\begin{equation*}\label{2.14-1}
	0<\underline{a}\leq e^{-\v p^{\v}_{in}}\leq \bar{a},~~\mbox{and}~~0<\underline{b}\leq e^{\t^\v_{in}}\leq \bar{b},
\end{equation*}
where $\underline{a}, \bar{a}, \underline{b}$ and $\bar{b}$ are given constants independent of $\v$. Then the uniform estimates for \eqref{2.1} are:

\begin{theorem}[Uniform Estimate for Ill-Prepared Data]\label{thm1.1}
	Let $s\geq 4$ be an integer, and assume that the initial data $(p^\v_{in},u^\v_{in},\t^\v_{in})$ satisfy 
	\begin{align}\label{initial}
		\|(p^\v_{in},u^\v_{in},\t^\v_{in}-\tilde{\t})\|^2_{H^{s,\v}}\leq \hat{C}_0<\infty, ~~\mbox{for}~~\v\in(0,1],
	\end{align}
	where $\hat{C}_0$ is  independent of $\v$.
	Then there exist positive constants $T_0$ and $\v_0$ depending only on $\hat{C}_0$ and $|\t_+-\t_-|$ such that the Cauchy problem \eqref{2.1}, \eqref{2.1-1} has a unique smooth solution $(p^\v,u^\v,\t^\v)$ satisfying 
	\begin{align}\label{ueid}
		\|(p^\v,u^\v,\t^\v-\tilde{\t})(t)\|^2_{H^{s,\v}} \leq \tilde{C}_0, ~~\mbox{for}~~t\in [0,T_0],~\v\in(0,\v_0],
	\end{align}
where the constant $\tilde{C}_0>0$ depends only  on $\hat{C}_0$ and $|\t_+-\t_-|$.
\end{theorem}

\begin{remark}
	The time derivatives of the initial data  in \eqref{initial} are defined through the system  \eqref{2.1}.
\end{remark}

\begin{remark}
	The wave strength $|\t_+-\t_-|$ (or equivalently $|\r_+-\r_-|$) is allowed to be large. But the constants $T_0$ and $\v_0$ may depend on the wave strength.
\end{remark}

Then we have the following low Mach limit.
\begin{theorem}[Low Mach Limit for Ill-Prepared Data]\label{ipC}
	Under the assumptions of Theorem \ref{thm1.1} and further assume that the initial data satisfy 
	\begin{equation}
		(p^\v_{in},u^\v_{in},\t^\v_{in}-\tilde{\theta})\rightarrow (p_{in},u_{in},\t_{in}-\tilde{\theta}),~~\mbox{in}~~H^s(\mathbb{R})~~ \mbox{as} ~~\v\rightarrow0,
	\end{equation}
	\begin{equation}\label{DCC}
		|\theta_{in}(x)-\theta_+|\leq C x^{-1-\s}, \mbox{ for } x\in [1, +\infty),
	\end{equation}
where $\s$ and $C$ are some positive constants. 
	Then the solutions $(p^\v,u^\v,\t^\v)$ of \eqref{2.1} converge strongly in $L^2(0, T_0; H^{s'}_{loc}(\mathbb{R}))$ for all $s'<s$ to  $(0, \overline{u},\overline{\t})$, where $( \overline{u},\overline{\t})$ is the unique solution  of \eqref{2.8} with the initial data $(w_{in},\t_{in})$, where $w_{in}$ is determined by 
	\begin{equation}
		w_{in}=\f12\k e^{\t_{in}}\partial_x\t_{in}.
	\end{equation}
	
\end{theorem}

\begin{remark}
In the proof of Theorem \ref{ipC}, we extend the convergence lemma of \cite{Alazard1, MS} to the different asymptotic states at infinity. 
\end{remark}

\begin{remark}
	If $\t_{in}^\v-\tilde{\t}\rightarrow 0$ as $\v\rightarrow0$, then one must have that $\bar{\t}=\tilde{\t}, \bar{u}=\f12\k e^{\tilde\t}\partial_x\tilde\t$ which is indeed the diffusive wave solution constructed in \eqref{1.15}.
\end{remark}

\begin{remark}
	The condition \eqref{DCC} can also be replaced by $|\theta_{in}(x)-\theta_-|\leq C |x|^{-1-\s}, \mbox{ for } x\in (-\infty,-1]$.
\end{remark}

The paper is organized as follows. In Section 2 we consider the low Mach limit for the well-prepared data and in Section 3  for the ill-prepared data.

\

\noindent\textbf{Notations:} Throughout this paper, the positive
generic constants independent of $\v$ are denoted by
$c$, $C$, $C_i$ ($i\in\mathbb{N}$). And we will use $\|\cdot\|$ to denote the standard
$L^2(\mathbb{R})$ norm, and $\|\cdot\|_{H^i}$ ($i\in\mathbb{N}$) to
denote the Sobolev $H^i(\mathbb{R})$ norm. Sometimes, we also use
$O(1)$ to denote a uniform bounded constant independent  of
$\v$.

\section{Well-prepared data}

This section is devoted to the low Mach limit for the well-prepared data. For simplicity, we omit the superscript $\v$ of the variables throughout this section.  We first reformulate the system \eqref{N-1.1} as follows. \\[1mm]

\noindent\underline{\bf Reformulation of the  system \eqref{N-1.1}:}\\[2mm]
 Set the scaling
\begin{equation}\label{N-2.2}
	y=\f{x}{\v},~~~\tau=\f{t}{\v^2}.
\end{equation}
In the following, we will also use the notations $(v,u,\mathcal{T})(\tau,y)$ and $(\tilde v,\tilde u,\tilde{\mathcal{T}})(\tau,y)$, etc., in the scaled independent variables. Set the perturbation around the profile $(\tilde v,\tilde u,\tilde{\mathcal{T}})(\tau,y)$ by
\begin{equation}\label{N-2.3}
(\phi,\psi,\omega,\zeta)(\tau,y)=\Big(v-\tilde v, \v u-\v\tilde u, \mathcal{T}+\f12|\v u|^2-\tilde{\mathcal{T}}-\f12|\v\tilde u|^2, \mathcal{T}-\tilde{\mathcal{T}}\Big)(\tau,y),
\end{equation}
and introduce 
\begin{eqnarray}\label{N-2.4}
	(\Phi,\Psi, \tilde W)(\tau,y)=\int_{-\infty}^{y}( \phi, \psi, \omega)(\tau,z)dz.
\end{eqnarray}
From  \eqref{N-1.7-1}, we note that
$(\Phi,\Psi,\tilde{W})(0,\pm\infty)=0$. And it follows from the equations of $(\Phi,\Psi,\tilde{W})$ in \eqref{N-2.5} below that $(\Phi,\Psi,\tilde{W})(\tau,\pm\infty)\equiv0$ for $\tau>0$. 

Subtracting \eqref{N-1.5} from \eqref{N-1.1} and integrating the resulting system yield the following integrated system for $(\Phi,\Psi,\tilde W)$:
\begin{equation}\label{N-2.5}
	\begin{cases}
		\Phi_{\tau}-\Psi_y=0,\\
		\Psi_\tau +P-\tilde{P}=\tilde\m\big(\f{1}{v}\v u_y-\f{1}{\tilde v}\v\tilde u_y\big)-\tilde R_1,\\
		\tilde W_\tau+ \v P u-\v\tilde P\tilde u=\k\big(\f{1}{v}\mathcal{T}_y-\f{1}{\tilde v}\tilde{\mathcal{T}}_y\big)+\tilde\mu\v^2\big(\f{1}{v}uu_y-\f{1}{\tilde v}\tilde u\tilde u_y\big)-\v\tilde R_2.
	\end{cases}
\end{equation}
We note that the term involving time derivative in  $\eqref{N-2.5}_3$ is $\tilde{W}_\tau$, while the dissipation term of $\eqref{N-2.5}_3$ is in the form $\k\big(\f{1}{v}\mathcal{T}_y-\f{1}{\tilde v}\tilde{\mathcal{T}}_y\big)$ which is closely related to the perturbation of temperature. To write the dissipation equation  $\eqref{N-2.5}_3$ in a unified way and make the linearized formulation more convenient to do the energy estimates,  instead of the variable $\tilde W$ which is the anti-derivative of the total energy, it is more convenient to introduce another variable $W$ related to the temperature, \textit{i.e.}, 
\begin{equation}\label{N-2.6}
	W=\tilde W-\v\tilde u\Psi,
\end{equation}
which yields that 
\begin{equation}\label{N-2.6-1}
	\zeta=W_y-Y,~~\mbox{with}~~Y=\f12|\Psi_y|^2-\v\tilde{u}_y\Psi.
\end{equation}
In terms of the new variable $W$, one linearize the system \eqref{N-2.5} as
\begin{equation}\label{N-2.7}
	\begin{cases}
		\Phi_{\tau}-\Psi_y=0,\\
		\Psi_\tau - \f{1}{\tilde v}\Phi_y +\f{1}{\tilde v}W_y=\tilde{\m}\f{1}{\tilde v} \Psi_{yy}+Q_1,\\
		W_\tau+ \Psi_y=\f{\k}{\tilde v}W_{yy}+Q_2,
	\end{cases}
\end{equation}
where
\begin{eqnarray}
	&& Q_1=\tilde{\m}\big(\f{1}{v}-\f{1}{\tilde v}\big)\v u_y+J_1+\f{1}{\tilde v}Y-\tilde R_1,\label{N-2.9}\\
	&& Q_2=\k\big(\f{1}{v}-\f{1}{\tilde v}\big)\mathcal{T}_y+\tilde{\m}\f{\v u_y}{v}\Psi_y+J_2-\v\tilde{u}_\tau\Psi-\f{\k}{\tilde v}Y_y-\v\tilde R_2+\v\tilde{u}\tilde R_1\label{N-2.10},
\end{eqnarray}
and 
\begin{equation*}
\begin{cases}
J_1=\f{\tilde{P}-1}{\tilde v}\Phi_y-[P-\tilde P+\f{\tilde{P}}{\tilde v}(v-\tilde v)-\f{1}{\tilde v}(\mathcal{T}-\tilde{\mathcal{T}})]\\
~~~=O(1)(|\Phi_y|^2+|W_y|^2+Y^2+|\v\tilde u|^4), \\[1mm]
J_2=(1-P)\Psi_y=O(1)(|(\Phi_y,\Psi_y,W_y)|^2+Y^2+|\v\tilde u|^4). 
\end{cases}
\end{equation*}

To prove Theorem \ref{Nthm1.1}, one needs only to show the following  {\it a priori} estimate in the scaled independent variables and employ the continuity argument since 
the local existence for compressible Navier-Stokes system is known. 
\begin{proposition}[A Priori Estimates]\label{Npro1.1}
	Assume that $(\Phi,\Psi, W)$ is a smooth solution of  \eqref{N-2.7} with zero initial data in the time interval $\tau\in [0,T]$. There exist constants $\d_1>0$ and $\v_1>0$ such that if $\d\leq \d_1$ and $\v\leq \v_1$, the following estimates hold
	\begin{equation*}
		\begin{cases}
			\|(\Phi,\Psi, W)(\tau)\|^2_{L^\infty_y}\leq C\sqrt{\d}\v,\\
			\|(\phi,\psi, \zeta)(\tau)\|^2_{L^2_y}\leq C\sqrt{\d}\v^2(1+\v^2\tau)^{-1+C_0\sqrt{\d}},\\
			\|(\phi_y,\psi_y, \zeta_y)(\tau)\|^2_{L^2_y}\leq C\sqrt{\d}\v^3(1+\v^2\tau)^{-\f32+C_0\sqrt{\d}}.
		\end{cases}
	\end{equation*}
\end{proposition}

To obtain the above {\it a priori estimates}, one assumes the following \textit{a priori} assumption:
\begin{equation}\label{N-2.13}
	N(T)=\sup_{0\leq \tau\leq T}\Big\{\f{1}{\v}\|(\Phi,\Psi, W)(\tau)\|^2_{L^\infty}+\f{1}{\v^2}\|(\phi,\psi, \zeta)(\tau)\|^2_{L^2}+\f{1}{\v^3}\|(\phi_y,\psi_y, \zeta_y)(\tau)\|^2_{L^2}\Big\}\leq \chi^2,
\end{equation}
where $\chi$ is a small constant determined later. It is noted that the local existence theorem yields that \eqref{N-2.13} holds for some positive time which may depend on $\v>0$ and $\d>0$.  We also point out that once we have proved the above Proposition \ref{Npro1.1}, then it implies that \eqref{N-2.13} is valid, and hence we justify the use of above {\it a priori} assumption \eqref{N-2.13}.

\

\noindent\underline{\bf Basic Estimates:}\\[2mm]
We start with the elementary energy estimates. Multiplying $\eqref{N-2.7}_1$ by $\Phi$, $\eqref{N-2.7}_2$ by $\tilde v\Psi$, $\eqref{N-2.7}_3$ by $W$ respectively,  and adding all the resulting equations, one can obtain that 
\begin{eqnarray}\label{N-2.14}
	&&\left(\f12\Phi^2+\f12\tilde{v}\Psi^2+\f{1}{2}W^2\right)_\tau+\tilde\mu\Psi_y^2+\f{\k}{\tilde v}W_y^2\nonumber\\
	&&=\f12\tilde{v}_\tau\Psi^2+\tilde{v}Q_1\Psi+Q_2W-\Big(\f{\k}{\tilde v}\Big)_yW_yW+(\cdots)_y.
\end{eqnarray}
Define
\begin{equation*}\label{N-2.15}
	m=(\Phi,\Psi, W)^t,
\end{equation*}
where $(\cdot,\cdot,\cdot)^t$ is the transpose of the vector $(\cdot,\cdot,\cdot)$, then from \eqref{N-2.7}, one has 
\begin{equation}\label{N-2.16}
	m_\tau+A_1m_y=A_2m_{yy}+A_3,
\end{equation}
where
\begin{equation*}\label{N-2.17}
	A_1=\left(\begin{array}{cccc}& 0 &-1&0\\&-\frac{1}{\tilde{v}}&0&\frac{1}{\tilde{v}}\\
		& 0 & 1 & 0\end{array}\right), \quad
	A_2=\left(\begin{array}{cccc}
		&0&0&0\\ & 0 &\frac{\tilde\mu}{\tilde{v}}& 0\\
		&0&0&\frac{\k}{\tilde{v}}\end{array}\right),~~
	A_3=(0,Q_1,Q_2)^t.
\end{equation*}
It follows from a direct computation that the eigenvalues of the matrix $A_1$
are $\lambda_1,0,\lambda_3$ with
$\lambda_3=-\lambda_1=\sqrt{\frac{2}{\tilde{v}}}$. The
corresponding normalized left and right eigenvectors can be chosen
as
\begin{eqnarray*}
	l_1=\f12(-1,-\frac{2}{\lambda_3},1),~
	l_2=\sqrt{\f{1}{2}}(1,0,1),~
	l_3=\f12(-1,\frac{2}{\lambda_3},1),\label{N-2.19}\\
	r_1=\f12(-1,-\lambda_3,1)^t,~
	r_2=\sqrt{\f{1}{2}}(1,0,1)^t,~r_3=\f12(-1,\lambda_3,1)^t \label{N-2.20}
\end{eqnarray*}
so that 
\begin{equation*}\label{N-2.21}
	l_ir_j=\delta_{ij},~i,j=1,2,3,\quad LA_1R=\Lambda=\left(\begin{array}{cccc}&\lambda_1&0&0\\&0&0&0\\
		&0&0&\lambda_3\end{array}\right),
\end{equation*}
where
$$L=(l_1,l_2,l_3)^t,~ R=(r_1,r_2,r_3).$$
Let
\begin{equation}\label{N-2.22}
	B=Lm=(b_1,b_2,b_3),
\end{equation}
then multiplying the
equations \eqref{N-2.16} by the matrix $L$ yields that
\begin{equation}\label{N-2.23}
	B_\tau+\Lambda B_y=LA_2RB_{yy}+2LA_2R_yB_y+[(L_\tau+\Lambda
	L_y)R+LA_2R_{yy}]B+LA_3.
\end{equation}
We shall use a weighted energy method to derive the
intrinsic dissipation. Without loss of generality, we assume that $\hat{\mathcal{T}}_{y}> 0$. The case that $\hat{\mathcal{T}}_y< 0$ can be treated in the same way. 
Let $\mathcal{T}_1=\f{\hat{\mathcal{T}}}{\mathcal{T}_+}$, then $|\mathcal{T}_1-1|\le
C\delta$. Multiplying \eqref{N-2.23} by
$\tilde{B}=(\mathcal{T}_1^Nb_1,b_2,\mathcal{T}_1^{-N}b_3)$ with a large positive integer
$N$ which will be chosen later, one can get that 
\begin{equation}\label{N-2.24}
	\begin{array}{ll}
		{\displaystyle\left(\frac{\mathcal{T}_1^N}2b_1^2+\frac12b_2^2+\frac{\mathcal{T}_1^{-N}}{2}b_3^2\right)_\tau
			-\left(\frac{\mathcal{T}_1^N}2\right)_\tau b_1^2-\left(\frac{\mathcal{T}_1^{-N}}2\right)_\tau b_3^2
			+\tilde{B}_yA_4B_y+\tilde{B}A_{4y}B_y}\\
		{\displaystyle-\frac{\mathcal{T}_1^{N-1}}2(N\lambda_1\mathcal{T}_{1y}+\mathcal{T}_1\lambda_{1y})b_1^2
			+\frac{\mathcal{T}_1^{-N-1}}2(N\lambda_3\mathcal{T}_{1y}-\mathcal{T}_1\lambda_{3y})b_3^2+(\cdots)_y}\\[0.3cm]
		={\displaystyle
			2\tilde{B}LA_2R_yB_y+\tilde{B}[L_\tau R+LA_2R_{yy}]B+\tilde{B}\Lambda
			L_yRB+\tilde{B}LA_3.}\end{array}
\end{equation}
A direct computation shows that the symmetric matrix $LA_2R=A_4$ is nonnegative. Define
\begin{align}
	E_1&=\int\Big(\f12\Phi^2+\f12\tilde{v}\Psi^2+\f{1}{2}W^2\Big)dy+\int\Big(\frac{\mathcal{T}_1^N}2b_1^2
	+\frac12b_2^2+\frac{\mathcal{T}_1^{-N}}{2}b_3^2\Big)dy,\label{N-2.25}\\
	K_1&=\int\Big(\tilde\mu\Psi_{y}^2
	+\frac{\k}{\tilde{v}}W_y^2+{B}_yA_4B_y\Big)dy.\label{N-2.26}
\end{align}
Note that
\begin{align}\label{N-2.27}
\left|\int(\tilde{B}-B)_yA_4B_ydy\right|&\leq 
			C\delta\int|B_y|^2dy+C\delta\int|\hat{\mathcal{T}}_y|^2|B|^2dy\nonumber\\
		&\leq C\v^2\delta(1+t)^{-1}E_1+C\delta
			K_1+C\delta\int|\Phi_y|^2dy.
\end{align}
Similarly, the terms  $\tilde{B}A_{4y}B_y$, $\tilde{B}LA_2R_yB_y$ and
$\tilde{B}[L_\tau R+LA_2R_{yy}]B$ in  \eqref{N-2.24} satisfy the same estimate. For
$\tilde{B}\Lambda L_yRB$ and $\tilde{B}LA_3$, we need to use the
explicit presentation. By the choice of the characteristic matrix
$L$ and $R$, one has that 
\begin{equation*}\label{N-2.28}
	\Lambda L_yR=\frac12\lambda_{3y}\left(\begin{array}{cccc}&1&0&-1\\&0&0&0\\
		&1&0&-1\end{array}\right),
	~~\mbox{and}~~
	LA_3=\left(\displaystyle\f12Q_2-\frac{Q_1}{\lambda_3},~ \sqrt{\f12}Q_2,~
	\f12Q_2+\frac{Q_1}{\lambda_3}\right)^t,
\end{equation*}
which yields immediately that
\begin{eqnarray}\label{N-2.30}
	|\tilde{B}\Lambda L_yRB|\leq C|\lambda_{3y}|(b_1^2+b_3^2),~~~\tilde{B}LA_3\leq C |(b_1,b_3)|\cdot|(Q_1,Q_2)|+C|b_2Q_2|.
\end{eqnarray}
Integrating \eqref{N-2.14} over $\mathbb{R}$ with respect to $y$,  using \eqref{N-2.24}-\eqref{N-2.27}, \eqref{N-2.30} and the Cauchy inequality, and choosing $N$ sufficiently large, we obtain that 
\begin{equation}\label{N-2.31}
	E_{1\tau}+\frac12	K_1+2\int|\hat{\mathcal{T}}_{y}|(b_1^2+b_3^2)|dy\le
	C\v^2\delta(1+t)^{-1}(E_1+1)+C\delta K_1
	+C\delta\int\Phi_y^2dy+I,
\end{equation}
with
\begin{equation}\label{N-2.32}
	I=C\int|(b_1,b_3)|\cdot|(Q_1,Q_2)|+|b_2Q_2|dy.
\end{equation}
Here, in the proof of \eqref{N-2.31}, we have used the fact
\begin{eqnarray}\label{N-2.33}
	\Psi=\f12 \l_3(b_3-b_1),
\end{eqnarray}
and for $N$ large enough that
\begin{equation}\label{N-2.34}
	-\frac12\mathcal{T}_1^{N-1}(N\lambda_1\hat{\mathcal{T}}_{1y}+2\mathcal{T}_1\lambda_{1y})b_1^2+\frac12\mathcal{T}_1^{-N-1}
	(N\lambda_3\mathcal{T}_{1y}-2\mathcal{T}_1\lambda_{3y})b^2_3-\tilde{B}\Lambda L_yRB
	\geq3|\hat{\mathcal{T}}_{y}|(b_1^2+b_3^2).
\end{equation}
Although $Q_1$ contains the term $\tilde{R}_1$ with the decay rate $\v^2(1+t)^{-1}$, the terms in \eqref{N-2.32} involving $Q_1$ can be estimated by the intrinsic dissipation on $b_1$ and $b_3$ as shown later. The terms involving $Q_2$ contains the term $\v\tilde{R}_2$ which has a better decay rate $\v^3(1+t)^{-\f32}$ and can be estimated directly. For brevity, we only estimate $\int|Q_1b_1|dy$ and $\int |Q_2b_2|dy$ as follows for
illustration.

\

\noindent\underline{Estimation on $\int|Q_1b_1|dy$:}

\

It follows  from \eqref{N-2.9} that
\begin{equation}\label{N-2.35}
	\int|Q_1b_1|dy\leq \int\Big|\Big[\Big(\f{\tilde{\m}}{v}-\f{\tilde{\m}}{\tilde v}\Big)u_y+J_1+\f{Y}{\tilde v}\Big]\cdot b_1\Big|+|\tilde R_1b_1|dy=:I_1+I_2.
\end{equation}
A direct calculation shows that 
\begin{align}
	\int\big|\Big(\f{\tilde{\m}}{v}-\f{\tilde{\m}}{\tilde v}\Big)\v u_y\cdot b_1\big|dy
	&\leq C (\chi+\d) (K_1+\|\Phi_y\|_{L^2}^2+\|\psi_y\|_{L^2}^2)+C\d\v^4(1+t)^{-2} \|\Phi\|_{L^2}^2\nonumber\\
	&\qquad+C\d\v^2(1+t)^{-1}E_1+C\d\v^5(1+t)^{-\f52},\nonumber
\end{align}
and
\begin{align}
\int\Big(|J_1|+|\f{Y}{\tilde v}|\Big)\cdot|b_1|dy
&\leq C(\d+\chi) (K_1+\|\Phi_y\|_{L^2}^2+\|\psi_y\|_{L^2}^2)+C\d\v^4(1+t)^{-2} \|\Phi\|_{L^2}^2\nonumber\\
&~~~~+C\d\v^2(1+t)^{-1}E_1+C\d\v^5(1+t)^{-\f52},\nonumber
\end{align}
which yields immediately that
\begin{equation}\label{N-2.38}
	I_1\leq C(\d+\chi) (K_1+\|\Phi_y\|_{L^2}^2+\|\psi_y\|_{L^2}^2)+  +C\d\v^2(1+t)^{-1}E_1+C\d\v^5(1+t)^{-\f52}.
\end{equation}
On the other hand, it follows from \eqref{N-1.6} and the Cauchy inequality that
\begin{equation*}\label{N-2.39}
	\int|\tilde{R}_1b_1|dy\leq \d\int|\hat{\mathcal{T}}_y|b_1^2dy+C\d\v^2(1+t)^{-1},
\end{equation*}
which, together with \eqref{N-2.38}, yields that
\begin{equation}\label{N-2.40}
	\int|Q_1b_1|dy\leq \d\int|\hat{\mathcal{T}}_y|b_1^2dy+C(\d+\chi) (K_1+\|\Phi_y\|_{L^2}^2+\|\psi_y\|_{L^2}^2)+C\d\v^2(1+t)^{-1}(E_1+1).
\end{equation}

\

\noindent\underline{Estimation on $\int|Q_2b_2|dy$:}\\[2mm]
It follows from \eqref{N-2.10} that
\begin{equation}\label{N-2.41}
	\int|Q_2b_2|dy\leq \int\big|\Big(\f{\k}{v}-\f{\k}{\tilde v}\Big)\mathcal{T}_y+\f{\tilde\mu}{v}\v u_y\Psi_y+J_2-\v\tilde{u}_\tau\Psi-\f{\k}{\tilde v}Y_y-\v\tilde R_2+\v\tilde{u}\tilde R_1\big|\cdot|b_2|dy.
\end{equation}
The Cauchy inequality yields that
\begin{align}\label{N-2.42}
	&\int\big|\Big(\f{\k}{v}-\f{\k}{\tilde v}\Big)\mathcal{T}_y+\f{\tilde\mu}{v}\v u_y\Psi_y\big|\cdot|b_2|dy\nonumber\\
	&\leq C(\d+\chi)(K_1+\|\Phi_y\|_{L^2}^2)+C\chi \|(\psi_y,\zeta_y)\|_{L^2}^2+C\d\v^2(1+t)^{-1}E_1,
\end{align}
and
\begin{align}\label{N-2.43}
	&\int\big|J_2-\v\tilde{u}_\tau\Psi-\f{\k}{\tilde v}Y_y-\v\tilde R_2+\v\tilde{u}\tilde R_1\big|\cdot|b_2|dy
	\nonumber\\
	&\leq C(\d+\chi)(K_1+\|\Phi_y\|_{L^2}^2)+C\chi \|\psi_y\|_{L^2}^2+C\d\v^2(1+t)^{-1}(E_1+1).
\end{align}
Thus combining \eqref{N-2.42}, \eqref{N-2.43} and using the Cauchy inequality, one has that
\begin{eqnarray}\label{N-2.44}
	\int|Q_2b_2|dy\leq C(\d+\chi) (K_1+\|\Phi_y\|_{L^2}^2)+C\chi \|(\psi_y,\zeta_y)\|_{L^2}^2+C\d\v^2(1+t)^{-1}(E_1+1).
\end{eqnarray}
Substituting  \eqref{N-2.40}, \eqref{N-2.44} into \eqref{N-2.31}, one obtains that 
\begin{align}\label{N-2.45}
	& E_{1\tau}+\frac14
	K_1+\int|\hat{\mathcal{T}}_{y}|(b_1^2+b_3^2)|dy\nonumber\\
	&\leq C\delta\v^2(1+t)^{-1}(E_1+1)+C(\d+\chi) (\|\Phi_y\|_{L^2}^2+\|(\psi_y,\zeta_y)\|_{L^2}^2).
\end{align}

\

Note that the norm  $\|\Phi_y\|_{L^2}^2$ is not included in $K_1$. To complete the lower-order inequality, we need to estimate $\Phi_y$. Using $\eqref{N-2.7}_1$, we rewrite  $\eqref{N-2.7}_2$ to be the following form 
\begin{eqnarray}\label{N-2.46}
	\frac{\tilde{\m}}{\tilde{v}}\Phi_{y\tau}-\Psi_{\tau}+\frac{1}{\tilde{v}}\Phi_y=\frac{1}{\tilde{v}}W_y-Q_1.
\end{eqnarray}
Multiplying \eqref{N-2.46} by $\Phi_y$ yields that
\begin{eqnarray}\label{N-2.47}
	\Big(\frac{\tilde\mu}{2\tilde{v}}\Phi_y^2\Big)_\tau-\Big(\frac{\tilde\mu}{2\tilde{v}}\Big)_\tau\Phi_y^2-\Phi_y\Psi_{\tau}+\frac{1}{\tilde{v}}\Phi_y^2=\Big(\frac{1}{\tilde{v}}W_y-Q_1\Big)\Phi_y.
\end{eqnarray}
Noting 
\begin{equation*}\label{N-2.48}
	\Psi_\tau\Phi_y=(\Psi\Phi_y)_\tau-(\Psi\Phi_\tau)_y+\Psi_y^2,
\end{equation*}
which, together with \eqref{N-2.47} and  integrating by parts, yields that 
\begin{equation}\label{N-2.49}
\Big(\int\frac{\tilde{\m}}{2\tilde{v}}\Phi_y^2-\Phi_y\Psi dy\Big)_\tau+\int\frac{1}{2\tilde{v}}\Phi_y^2dy
\leq C\int\Psi_{y}^2+W_y^2 dy+C\delta(1+t)^{-3/2}+\int Q_1^2dy.
\end{equation}
On the other hand, it follows from \eqref{N-2.9} and \eqref{N-2.13}  that
\begin{align}\label{N-2.50}
	\int Q_1^2dy\leq C(\d+\chi)(K_1+\|\Phi_y\|^2+\|\psi_y\|^2)+C\delta\v^2(1+t)^{-1}E_1+C\delta\v^2(1+t)^{-1}.
\end{align}
Substituting \eqref{N-2.50} into \eqref{N-2.49}, one has that 
\begin{equation}\label{N-2.51}
	\Big(\int\frac{\tilde\mu}{2\tilde{v}}\Phi_y^2-\Phi_y\Psi dy\Big)_\tau+\int\frac{1}{4\tilde{v}}\Phi_y^2dy\leq C_1K_1+C_1\delta\v^2(1+t)^{-1}(E_1+1)+C_1(\d+\chi)\|\psi_y\|^2.
\end{equation}
By choosing $\tilde C_1$ large enough, it holds that 
\begin{equation}\label{N-2.52}
	\tilde C_1E_1+\int\frac{\tilde\mu}{2\tilde{v}}\Phi_y^2-\Phi_y\Psi dy\geq \f12 \tilde C_1E_1+\int\frac{\tilde\mu}{4\tilde{v}}\Phi_y^2dy,~~\mbox{and}~~ \f14\tilde C_1-C_1\geq \f18\tilde C_1.
\end{equation}
Thus, multiplying \eqref{N-2.45} by $\tilde{C}_1$ and using \eqref{N-2.52}, one obtains the low order estimates:
\begin{lemma}\label{Nthm2.1}
	If $\d$ and $\chi$ are suitably small, then it holds that 
	\begin{equation*}\label{N-2.53}
		E_{2\tau}+K_2+\int|\hat{\mathcal{T}}_{y}|(b_1^2+b_3^2)|dy\leq C\delta\v^2(1+t)^{-1}(E_1+1)+C(\d+\chi) \|(\psi_y,\zeta_y)\|^2,
	\end{equation*}
	where
	\begin{equation*}\label{N-2.54}
		E_2=\tilde{C}_1E_1+\int\frac{\tilde\mu}{2\tilde{v}}\Phi_y^2-\Phi_y\Psi dy,~~
		K_2=\f18\tilde C_1K_1+\int\frac{1}{8\tilde{v}}\Phi_y^2dy.
	\end{equation*}
\end{lemma}

\

\noindent\underline{\bf Derivative Estimates:}\\[2mm]
To obtain the estimate for the first-order derivative of $\|(\Phi_y,\Psi_y,W_y)\|$, we shall
use an energy estimate based on the convex entropy of the compressible Navier-Stokes equations. Applying $\partial_y$ to \eqref{N-2.5} yields that 
\begin{equation}\label{N-2.55}
	\begin{cases}
		\phi_{\tau}-\psi_y=0,\\
		\psi_\tau +(P-\tilde{P})_y=\v(\f{\tilde\mu}{v}u_y-\f{\tilde\mu}{\tilde v}\tilde{u}_y)_y -\tilde{R}_{1y},\\
		\zeta_\tau+ \v(Pu_y-\tilde{P}\tilde{u}_y)=(\f{\k}{v}\mathcal{T}_y-\f{\k}{\tilde v}\tilde{\mathcal{T}}_y)_y+Q_3,
	\end{cases}
\end{equation}
where
\begin{equation*}\label{N-2.56}
	Q_3=\f{\tilde\mu}{ v}\v^2u_y^2-\v^2\Big(\f{\tilde\mu}{\tilde v}\tilde{u}\tilde{u}_y\Big)_y+\v\tilde{P}_y\tilde{u}+\f12(\v^2\tilde{u}^2)_\tau-\v\tilde{R}_{2y}.
\end{equation*}

\begin{lemma}\label{Nthm2.2}
	If $\d$ and $\chi$ are suitably small, then it holds that 
	\begin{eqnarray}\label{N-2.71}
		E_{3\tau}+\f12K_3\leq C\d\v^2(1+t)^{-1}E_3+\v^2(1+t)^{-1}\int|\hat{\mathcal{T}}_y|(b_1^2+b_3^2)dy+C\d\v^4(1+t)^{-2},
	\end{eqnarray}
	where
	\begin{equation*}
		E_{3}=\int\f12\psi^2+\tilde{\mathcal{T}}\hat{\Phi}\Big(\frac{v}{\tilde{v}}\Big)
		+\tilde{\mathcal{T}}\hat{\Phi}\Big(\frac{\mathcal{T}}{\tilde{\mathcal{T}}}\Big)dy,
		~~~K_3=\int\frac{\tilde\mu}{v}\psi_y^2+\f34\f{\k}{v\mathcal{T}}\zeta_y^2dy.\label{N-2.73}
	\end{equation*}
	
\end{lemma}
\noindent\textbf{Proof}.
Multiplying $\eqref{N-2.55}_2$ by $\psi$ yields that 
\begin{equation*}\label{N-2.57}
	\Big(\f12\psi^2\Big)_\tau-(P-\tilde{P})\psi_{y}+\v\Big(\frac{\tilde\mu}{v}u_{y}
	-\frac{\mu(\tilde{\theta})}{\tilde{v}}\tilde{u}_{y}\Big)\psi_{y}=-\tilde{R}_{1y}\psi+(\cdots)_y.
\end{equation*}
Since $P-\tilde{P}=\tilde{\mathcal{T}}(\frac1v-\frac1{\tilde
	v})+\frac{\zeta}{v}$, as in \cite{Huang-2005}, it holds that 
\begin{align}\label{N-2.58}
	\Big(\f12\psi^2\Big)_\tau-\tilde{\mathcal{T}}\Big(\frac{1}{v}-\frac{1}{\tilde{v}}\Big)\phi_\tau-\frac{\zeta}{v}\psi_y+\frac{\tilde\mu}{v}\psi_y^2+\Big(\frac{\tilde\mu}{v}
	-\frac{\tilde{\m}}{\tilde{v}}\Big)\v\tilde{u}_{y}\psi_{y}=-\tilde{R}_{1y}\psi+(\cdots)_y.
\end{align}
Denote
\begin{equation}\label{N-2.59}
	\hat{\Phi}(s)=s-1-\ln s,
\end{equation}
then, it is straightforward to check that $\hat{\Phi}'(1)=0$ and $\hat{\Phi}(s)$ is
strictly convex around $s=1$. Moreover, it holds that
\begin{align}\label{N-2.60}
	\left(\tilde{\mathcal{T}}\hat{\Phi}\Big(\frac{v}{\tilde{v}}\Big)\right)_{\tau}&=\tilde{\mathcal{T}}_{\tau}\hat{\Phi}\Big(\frac{v}{\tilde{v}}\Big)
	+\tilde{\mathcal{T}}\Big(-\frac{1}{v}+\frac{1}{\tilde{v}}\Big)\phi_{\tau}+\tilde{\mathcal{T}}\Big(-\frac{v}{\tilde{v}^2}+\frac{1}{ \tilde{v}}\Big) \tilde{v}_{\tau}
	+\tilde{\mathcal{T}}\Big(-\frac{1}{v}+\frac{1}{\tilde{v}}\Big)\tilde{v}_{\tau}\nonumber\\
	&=\tilde{\mathcal{T}}\Big(-\frac{1}{v}+\frac{1}{\tilde{v}}\Big)\phi_{\tau}-\tilde{P}\hat{\Psi}\Big(\frac{v}{\tilde{v}}\Big)\tilde{v}_{\tau}
	+\tilde{v}\tilde{P}_{\tau}\hat{\Phi}\Big(\frac{v}{\tilde{v}}\Big),
\end{align}
with $\hat{\Psi}(s)=\hat{\Phi}(s^{-1})$. Therefore, it follows from \eqref{N-2.58} and  \eqref{N-2.60} that
\begin{eqnarray}\label{N-2.62}
	&&\Big(\f12\psi^2+\tilde{\mathcal{T}}\hat{\Phi}\Big(\frac{v}{\tilde{v}}\Big)\Big)_\tau
	+\tilde{P}\tilde{v}_{\tau}\hat{\Psi}\Big(\frac{v}{\tilde{v}}\Big)-\frac{\zeta}{v}\psi_y
	+\frac{\tilde\mu}{v}\psi_y^2+\Big(\frac{\tilde\mu}{v}
	-\frac{\tilde\mu}{\tilde{v}}\Big)\v\tilde{u}_{y}\psi_{y}\nonumber\\
	&&=\tilde{v}\tilde{P}_{\tau}\hat{\Phi}\Big(\frac{v}{\tilde{v}}\Big)-\tilde{R}_{1y}\psi+(\cdots)_y.
\end{eqnarray}
On the other hand, multiplying $\eqref{N-2.55}_3$ by $\f{\zeta}{\mathcal{T}}$ yields that
\begin{eqnarray}\label{N-2.63}
	\f{\zeta}{\mathcal{T}}\zeta_\tau+\v(Pu_y-\tilde{P}\tilde{u}_y)\f{\zeta}{\mathcal{T}}
	=\Big(\f{\k}{v}\mathcal{T}_y-\f{\k}{\tilde v}\tilde{\mathcal{T}}_y\Big)_y\f{\zeta}{\mathcal{T}}+Q_3\f{\zeta}{\mathcal{T}}.
\end{eqnarray}
A direct calculation gives that 
\begin{equation}\label{N-2.64}
	\f{\zeta}{\mathcal{T}}\zeta_\tau=\Big(\tilde{\mathcal{T}}\hat{\Phi}(\frac{\mathcal{T}}{\tilde{\mathcal{T}}})\Big)_{\tau}
	+O(1)|\tilde{\mathcal{T}}_\tau|\zeta^2,
\end{equation}
\begin{equation}\label{N-2.65}
	\v(Pu_y-\tilde{P}\tilde{u}_y)\f{\zeta}{\mathcal{T}}=\frac{\zeta}{v}\psi_y+(P-\tilde P)\v\tilde{u}_y\f{\zeta}{\mathcal{T}},
\end{equation}
and
\begin{equation}\label{N-2.66}
	\Big(\f{\k}{v}\mathcal{T}_y-\f{\k}{\tilde v}\tilde{\mathcal{T}}_y\Big)_y\f{\zeta}{\mathcal{T}}\leq -\f34\f{\k}{v\mathcal{T}}\zeta_y^2+O(1)|\tilde{\mathcal{T}}_y|^2|(\phi,\zeta)|^2
	+(\cdots)_y.
\end{equation}
Substitute \eqref{N-2.64}-\eqref{N-2.66} into \eqref{N-2.63}, one obtains that 
\begin{equation}\label{N-2.67}
	\Big(\tilde{\mathcal{T}}\hat{\Phi}\Big(\frac{\mathcal{T}}{\tilde{\mathcal{T}}}\Big)\Big)_{\tau}+\frac{\zeta}{v}\psi_y+\f34\f{\k}{v\mathcal{T}}\zeta_y^2
	\leq C\d\v^2(1+t)^{-1}|(\phi,\zeta)|^2+Q_3\f{\zeta}{\mathcal{T}}+(\cdots)_y.
\end{equation}
It follows from  \eqref{N-2.67} and \eqref{N-2.62} that
\begin{align}\label{N-2.68}
	&\Big(\f12\psi^2+\tilde{\mathcal{T}}\hat{\Phi}\Big(\frac{v}{\tilde{v}}\Big)
	+\tilde{\mathcal{T}}\hat{\Phi}\Big(\frac{\mathcal{T}}{\tilde{\mathcal{T}}}\Big)\Big)_\tau
	+\frac{\tilde\mu}{v}\psi_y^2+\f34\f{\k}{v\mathcal{T}}\zeta_y^2\nonumber\\
	&\leq C\d\v^2(1+t)^{-1}|(\phi,\zeta)|^2-\tilde{R}_{1y}\psi+Q_3\f{\zeta}{\mathcal{T}}+(\cdots)_y.
\end{align}
Using \eqref{N-2.33}, one has that 
\begin{equation}\label{N-2.69}
	\Big|\int\tilde{R}_{1y}\psi dy\Big|=\Big|\int\tilde{R}_{1yy}\Psi dy\Big|\leq \v^2(1+t)^{-1}\int|\hat{\mathcal{T}}_y|(b_1^2+b_3^2)dy+C\d\v^4(1+t)^{-2},
\end{equation}
and
\begin{align}\label{N-2.70}
	&\Big|\int Q_3\f{\zeta}{\t} dy\Big|\leq C\int|\zeta|\cdot|\psi_y|^2dy+C\d\v^4\int|\zeta|(1+t)^{-2}e^{-\f{cx^2}{1+t}}dy\nonumber\\
	&\leq C\chi\|\psi_y\|^2+C\d\v^2(1+t)^{-1}\|\zeta\|^2+C\d\v^5(1+t)^{-\f52}.
\end{align}
Integrating \eqref{N-2.68} with respect to $y$,  using \eqref{N-2.69} and \eqref{N-2.70}, then one obtains \eqref{N-2.71}. Thus, the proof of Lemma \ref{Nthm2.2} is completed. 
$\hfill\Box$

\

Since the norm  $\|\phi_y\|^2$ is not included in $K_3$. To complete the first derivative inequality, we follow the same way  for  $\Phi_y$. We rewrite the equation $\eqref{N-2.55}_2$ as
\begin{equation}\label{N-2.74}
	\frac{\tilde\mu}{\tilde{v}}\phi_{y\tau}-\psi_{\tau}+(P-\tilde P)_y=-\Big(\frac{\tilde\mu}{\tilde{v}}\Big)_y\phi_{y}
	-\Big(\Big(\f{\tilde\mu}{v}-\f{\tilde\mu}{\tilde{v}}\Big)\v u_y\Big)_y+\tilde{R}_{1y}.
\end{equation}
\begin{lemma}\label{Nthm2.3}
	If $\d$ and $\chi$ are suitably small, then it holds that 
	\begin{align}\label{N-2.78}
		&\Big(\int\frac{\tilde\mu}{2\tilde{v}}\phi_y^2-\phi_y\psi dy\Big)_\tau+\int
		\frac{\tilde{P}}{2\tilde{v}}\phi_y^2dy\nonumber\\
		&\leq C_2K_3+C_2\d\v^2(1+t)^{-1}E_3+C_2\d\v^5(1+t)^{-\frac52}
		+C_2(\|(\phi,\psi)\|^2_{L^\infty}+\|\psi_y\|^2)\|\psi_{yy}\|^2.
	\end{align}
\end{lemma}
\noindent\textbf{Proof}.
Multiplying \eqref{N-2.74} by $\phi_y$ yields that
\begin{equation*}\label{N-2.75}
\Big(\frac{\tilde\mu}{2\tilde{v}}\phi_{y}^2\Big)_\tau-\psi_{\tau}\phi_y-(P-\tilde{P})_y\phi_y 
=\Big(\frac{\tilde\mu}{2\tilde{v}}\Big)_\tau\phi_{y}^2+\Big\{-(\frac{\tilde\mu}{\tilde{v}})_y\psi_{y}
-[(\frac{\tilde\mu}{v}-\frac{\tilde\mu}{\tilde{v}})\v u_{y}]_y\Big\}\phi_y+\tilde{R}_{1y}\phi_y.
\end{equation*}
A direct calculation yields that 
\begin{equation*}\label{N-2.76}
-(P-\tilde{P})_y=\frac{\tilde{P}}{\tilde{v}}\phi_y-\frac{1}{\tilde{v}}\zeta_y+\Big(\frac{P}{v}-\frac{\tilde{P}}{\tilde{v}}\Big)v_y-\Big(\frac{1}{v}-\frac{1}{\tilde{v}}\Big)\mathcal{T}_y,
\end{equation*}
and
\begin{equation*}\label{N-2.77}
	\phi_y\psi_{\tau}=(\phi_y\psi)_\tau-(\phi_\tau\psi)_y+\psi_{y}^2.
\end{equation*}
Then combining the above three equations, one obtains that 
\begin{align}\label{N-2.78-1}
	&\Big(\int\frac{\tilde\mu}{2\tilde{v}}\phi_y^2-\phi_y\psi dy\Big)_\tau+\int
	\frac{\tilde{P}}{2\tilde{v}}\phi_y^2dy\nonumber\\
	&\leq C_2K_3+C_2\d\v^2(1+t)^{-1}E_3+C_2\d\v^5(1+t)^{-\frac52}
	+C_2(\|(\phi,\psi)\|^2_{L^\infty}+\|\psi_y\|^2)\|\psi_{yy}\|^2,
\end{align}
where we have used
\begin{equation*}\label{N-2.79}
	-\int(P-\tilde P)_y\phi_ydy\geq \int\f{3\tilde{P}}{4\tilde{v}}\phi_y^2dy-C\|\zeta_y\|^2-C\d\v^2(1+t)^{-1}E_3,
\end{equation*}
\begin{align}
	&\Big|\int(\frac{\tilde\mu}{\tilde{v}})_y\psi_{y}dy\Big|+
	\Big|\int[(\frac{\tilde\mu}{v}-\frac{\tilde\mu}{\tilde{v}})\v u_{y}]_y\phi_ydy\Big|
	+\Big|\int\tilde{R}_{1y}\phi_{y}dy\Big|\nonumber\\
	&\leq \int\f{\tilde{P}}{8\tilde{v}}\phi_y^2dy+CK_3+C\d\v^2(1+t)^{-1}E_3+C_2\d\v^5(1+t)^{-\frac52}\nonumber\\
	&~~~+\int|\phi_y|^2|\psi_y|+|\phi_y\psi_y\zeta_y|dy,\nonumber
\end{align}
and
\begin{align}
\int|\phi_y|^2|\psi_y|+|\phi_y\psi_y\zeta_y|dy&\leq C\|\phi_y\|^2\|\psi_y\|^{\f12}\|\psi_{yy}\|^{\f12}+C\|\phi_y\|\|\zeta_y\|\|\psi_y\|^{\f12}\|\psi_{yy}\|^{\f12}\nonumber\\
&\leq \int\f{\tilde{P}}{16\tilde{v}}\phi_y^2dy+C\|\zeta_y\|^2+C\|\psi_y\|^2\|\psi_{yy}\|^2.\nonumber
\end{align}
Therefore the proof of Lemma \ref{Nthm2.3} is completed. 
$\hfill\Box$

\

We now derive the higher order estimates. Applying $\partial_y$ to \eqref{N-2.55} yields that 
\begin{equation}\label{N-2.82}
	\begin{cases}
		\phi_{y\tau}-\psi_{yy}=0,\\
		\psi_{y\tau} +\f{1}{\tilde{v}}\zeta_{yy}-\f{1}{\tilde{v}}\phi_{yy}=\tilde\mu(\f{1}{v}\v u_y-\f{1}{\tilde v}\v\tilde{u}_y)_{yy}+Q_4 -\tilde{R}_{1yy},\\
		\zeta_{y\tau}+ \psi_{yy}=(\f{\k}{v}\mathcal{T}_y-\f{\k}{\tilde v}\tilde{\mathcal{T}}_y)_{yy}+Q_5,
	\end{cases}
\end{equation}
where
\begin{align}
	& Q_4=\f{\tilde{P}-1}{\tilde v}\phi_{yy}+\Big(\f{P}{v}-\f{\tilde{P}}{\tilde v}\Big)\phi_{yy}
	-2\Big(\f{1}{v^2}\mathcal{T}_yv_y-\f{1}{\tilde{v}^2}\tilde{\mathcal{T}}_y\tilde{v}_y\Big)
	+2\Big(\f{1}{v^3}v_y^2-\f{1}{\tilde{v}^3}\tilde{v}_y^2\Big)\nonumber\\
	&~~~~~+O(1)\Big(|\tilde{v}_{yy}||(\phi,\psi)|+|\phi\zeta_{yy}|\Big),\label{N-2.83}\\
	& Q_5=-\v\tilde{u}_{yy}(P-\tilde P)-\v(P_yu_y-\tilde{P}_y\tilde{u}_y)+(1-\tilde{P})\psi_{yy}+Q_{3y}.\label{N-2.84}
\end{align}

\begin{lemma}\label{Nlem2.4}
	If $\d$ and $\chi$ are suitably small, then it holds that 
	\begin{align}\label{N-2.95}
		E_{4\tau}+\f12K_4&\leq C\d\v(1+t)^{-\f12}\|(\phi_y,\psi_y,\zeta_y)\|^2+C\d\v^3(1+t)^{-\f32}E_3+C\|(\phi_y,\psi_y,\zeta_y)\|^{\f{10}{3}}\nonumber\\
		&~~~+C\|(\phi,\psi)\|\cdot\|(\phi_y,\psi_y,\zeta_y)\|^{3}
		+C\d\v^6(1+t)^{-3},
	\end{align}
	where
	\begin{equation*}\label{N-2.96}
		E_{4}=\int\f12\phi_y^2+\f12\tilde{v}\psi_y^2+\f12\zeta_y^2dy,
		~~K_4=\int\f{\tilde\mu\tilde{v}}{v}\psi_{yy}^2+\f{\k}{v}\zeta_{yy}^2dy.
	\end{equation*}
\end{lemma}
\noindent\textbf{Proof}. Multiplying $\eqref{N-2.82}_1$ by $\phi_y$, $\eqref{N-2.82}_2$ by $\tilde{v}\psi_y$, $\eqref{N-2.82}_3$ by $\zeta_y$ and adding the resulting equations, one obtains that 
\begin{align}\label{N-2.85}
	&\left(\f12\phi_y^2+\f12\tilde{v}\psi_y^2+\f12\zeta_y^2\right)_\tau
	+\f{\tilde\mu\tilde{v}}{v}\psi_{yy}^2+\f{\k}{v}\zeta_{yy}^2\nonumber\\
	&=\f12\tilde{v}_\tau\psi_y^2+\tilde{v}Q_4\psi_y -\tilde{R}_{1yy}\tilde{v}\psi_y+Q_5\zeta_y+J_3+J_4+(\cdots)_y,
\end{align}
where
\begin{align}
	& J_3=-\tilde{\m}\Big(\f{1}{v}\v u_y-\f{1}{\tilde v}\v\tilde{u}_y\Big)_{y}\tilde{v}_{y}\psi_y
	-\Big(\f{\tilde\mu}{v}\Big)_y\psi_{yy}\psi_y-\Big(\Big(\f{\tilde\mu}{v}-\f{\tilde\mu}{\tilde v}\Big)\v\tilde{u}_y\Big)_{y}\tilde{v}\psi_{yy},\nonumber\\
	& J_4=-\Big(\f{\k}{v}\Big)_y\zeta_{yy}\zeta_y-\Big(\Big(\f{\k}{v}-\f{\k}{\tilde v}\Big)\tilde{\t}_y\Big)_{y}\zeta_{yy}.\nonumber
\end{align}
It is direct to obtain
\begin{eqnarray}\label{N-2.88}
	&&\Big|\int\Big\{\f{\tilde{P}-1}{\tilde v}\phi_{yy}+\Big(\f{P}{v}-\f{\tilde{P}}{\tilde v}\Big)\phi_{yy}\Big\}\tilde{v}\psi_ydy\Big|
	=\Big|\int\Big\{(\tilde{P}-1)\psi_y+\tilde{v}\psi_y\Big(\f{P}{v}-\f{\tilde{P}}{\tilde v}\Big)\Big\}_y\phi_{y}dy\Big|\nonumber\\
	&&\leq(\d+\chi)\|\psi_{yy}\|^2+C\d\v^3(1+t)^{-\f32}\|(\phi_y,\psi_y,\zeta_y)\|^2
	+C\d\v^4(1+t)^{-2}E_3\nonumber\\
	&&~~~~~~~~~+C\|(\phi_y,\psi_y,\zeta_y)\|^{\f{10}{3}}+C\|(\phi,\psi)\|\|(\phi_y,\psi_y,\zeta_y)\|^{3},
\end{eqnarray}
\begin{align}\label{N-2.89}
	& 2\Big|\int\Big(\f{1}{v^2}\mathcal{T}_yv_y-\f{1}{\tilde{v}^2}\tilde{\mathcal{T}}_y\tilde{v}_y\Big)\tilde{v}\psi_ydy\Big|
	+2\Big|\int\Big(\f{1}{v^3}v_y^2-\f{1}{\tilde{v}^3}\tilde{v}_y^2\Big)\tilde{v}\psi_ydy\Big|\nonumber\\
	& \leq \chi\|\psi_{yy}\|^2+C\d\v(1+t)^{-\f12}\|(\phi_y,\psi_y,\zeta_y)\|^2
	+C\d\v^3(1+t)^{-\f32}E_3+C\|(\phi_y,\psi_y,\zeta_y)\|^{\f{10}{3}},
\end{align}
and
\begin{align}\label{N-2.90}
	\int\big(|\tilde{v}_{yy}||(\phi,\psi)|+|\phi\zeta_{yy}|\big)|\tilde{v}\psi_y|dy
	&\leq \chi\|\zeta_{yy}\|^2+C\d\v(1+t)^{-\f12}\|(\phi_y,\psi_y,\zeta_y)\|^2\nonumber\\
	&~~~+C\d\v^3(1+t)^{-\f32}E_3+C\|(\phi,\psi)\|\|(\phi_y,\psi_y,\zeta_y)\|^{3}.
\end{align}
Then it follows from \eqref{N-2.83} and \eqref{N-2.88}-\eqref{N-2.90} that
\begin{align}\label{N-2.91}
	&\Big|\int\tilde{v} Q_4\psi_ydy\Big|
	\leq (\d+\chi)\|(\psi_{yy},\zeta_{yy})\|^2+C\d\v(1+t)^{-\f12}\|(\phi_y,\psi_y,\zeta_y)\|^2+C\d\v^3(1+t)^{-\f32}E_3\nonumber\\
	&~~~~~~~~~~~~~~~~~~~~~~~~+C\|(\phi,\psi)\|\|(\phi_y,\psi_y,\zeta_y)\|^{3}+C\|(\phi_y,\psi_y,\zeta_y)\|^{\f{10}{3}}.
\end{align}
Using Cauchy inequality, it is easy to get that 
\begin{equation}\label{N-2.92}
\Big|\int \tilde{v}\tilde{R}_{1yy}\psi_ydy\Big|\leq C\d\v(1+t)^{-\f12}\|\psi_y\|^2+C\d\v^6(1+t)^{-3}.
\end{equation}
A direct calculation also yields that 
\begin{align}\label{N-2.93}
\Big|\int Q_5\zeta_ydy\Big|&\leq (\d+\chi)\|(\psi_{yy},\zeta_{yy})\|^2+C\d\v(1+t)^{-\f12}\|(\phi_y,\psi_y,\zeta_y)\|^2\nonumber\\
	&~~~+C\d\v^3(1+t)^{-\f32}E_3+C\|(\phi_y,\psi_y,\zeta_y)\|^{\f{10}{3}},
\end{align}
and
\begin{align}\label{N-2.94}
\int |J_3|+|J_4|dy&\leq (\d+\chi)\|(\psi_{yy},\zeta_{yy})\|^2+C\d\v(1+t)^{-\f12}\|(\phi_y,\psi_y,\zeta_y)\|^2\nonumber\\
&~~~~+C\d\v^3(1+t)^{-\f32}E_3+C\|(\phi_y,\psi_y,\zeta_y)\|^{\f{10}{3}}.
\end{align}
Integrating \eqref{N-2.85} with respect to $y$ and using \eqref{N-2.91}-\eqref{N-2.94}, one obtains \eqref{N-2.95}. Therefore the proof of Lemma \ref{Nlem2.4} is completed. $\hfill\Box$

\

We choose a large constants $\tilde{C}_2>1$ so that
\begin{equation}\label{N-2.97}
	\tilde{C}_2E_3+\int\f{2\mu(\tilde\t)}{3\tilde{v}}\phi_y^2-\psi\phi_ydy\geq \f12\tilde{C}_2E_3+\int\f{2\mu(\tilde\t)}{3\tilde{v}}\phi_y^2dy,~~\f12\tilde{C}_2-C_2\geq \f14\tilde{C}_2,
\end{equation}
and define
\begin{align}\label{N-2.98}
	E_{5}=\tilde{C}_2E_3+\int\f{2\mu(\tilde\t)}{3\tilde{v}}\phi_y^2-\psi\phi_ydy+\v^{-1}E_4,
	~~K_5=\f14\tilde{C}_2K_3+\int\f{\tilde{p}}{4\tilde{v}}\phi_{y}^2dy+\f12\v^{-1}K_4.
\end{align}
Thus, combining  \eqref{N-2.71}, \eqref{N-2.78} and  \eqref{N-2.95}, one obtains 
\begin{lemma}
	If $\d$ and $\chi$ are suitably small, then it holds that 
	\begin{align}\label{N-2.99}
		& E_{5\tau}+K_5\leq C_3\v^2(1+t)^{-1}\int|\hat{\mathcal{T}}_y|(b_1^2+b_2^2)dy+ C_3\d\v^2(1+t)^{-\f32}E_5+C_3\d\v^4(1+t)^{-2}.
	\end{align}
\end{lemma}

\

\

\noindent{\bf Proof of Proposition \ref{Npro1.1}:}
Let $\tilde{C}_3=1+C_3$, then from \eqref{N-2.51} and \eqref{N-2.99}, one gets that
\begin{equation*}\label{N-2.100}
	E_{6\tau}+K_6\leq C_0\sqrt{\d}\v^2(1+\v^2\tau)^{-1}(E_6+\sqrt\d).
\end{equation*}
where
\begin{equation*}\label{N-2.100-1}
	E_{6}=\tilde{C}_3E_2+\v^{-2}E_5,
	~~K_6=\f14\tilde{C}_3K_2+\f12\v^{-2}K_5.
\end{equation*}
Then the Gronwall's inequality yields that
\begin{equation}\label{N-2.101}
	E_{6}(\tau)\leq C_0\sqrt\d(1+\v^2\tau)^{C_0\sqrt\d},
	~~\int_{0}^{\tau} K_6(s)ds\leq C C_0\sqrt\d(1+\v^2\tau)^{C_0\sqrt\d}.
\end{equation}
Since $c_3\|(\Phi,\Psi,W)\|^2\leq E_6$ for some positive constant $c_3$, one obtains that
\begin{equation}\label{N-2.102}
	\|(\Phi,\Psi,W)\|^2 \leq C C_0\sqrt\d(1+\v^2\tau)^{C_0\sqrt\d}.
\end{equation}
Multiplying \eqref{N-2.99} by $(1+\v^2\tau)$, one has that 
\begin{eqnarray}\label{N-2.103}
	&&[(1+\v^2\tau)E_{5}]_\tau+(1+\v^2\tau)K_5\leq C\v^2\int|\hat{\mathcal{T}}_y|(b_1^2+b_2^2)dy+ C\v^2E_5+C_3\d\v^4(1+\v^2\tau)^{-1}\nonumber\\
	&&~~~~~~~~~~~~~~~~~~~~~~~~~~~~~~~~~~~~~~\leq C\v^2K_6+C_3\d\v^4(1+\v^2\tau)^{-1},
\end{eqnarray}
where we have used the fact that
\begin{align}
	E_{5}(\tau)&\leq C\|(\phi,\psi,\zeta)\|^2+C\v^{-1}\|(\phi_y,\psi_y,\zeta_y)\|^2\nonumber\\
	&\leq C\|(\Phi_y,\Psi_y,W_y)\|^2+C\v^{-1}\|(\phi_y,\psi_y,\zeta_y)\|^2+C\d\v^4(1+\v^2\tau)^{-\f32}\nonumber\\
	&\leq CK_6+C\d\v^4(1+\v^2\tau)^{-\f32}.\nonumber
\end{align}
Integrating \eqref{N-2.103} over $[0,\tau]$ and using \eqref{N-2.101}, one obtains that
\begin{equation}\label{N-2.105}
	E_{5}\leq C C_0\sqrt\d\v^2(1+\v^2\tau)^{-1+C_0\sqrt\d},~~~\int_{0}^{\tau}(1+\v^2s)K_5ds\leq C C_0\sqrt\d\v^2(1+\v^2\tau)^{C_0\sqrt\d}.
\end{equation}
Furthermore, it holds that
\begin{align}\label{N-2.106}
	E_{5}(\tau)&\geq c_4\|(\phi,\psi,\zeta)\|^2+c_4\v^{-1}\|(\phi_y,\psi_y,\zeta_y)\|^2\nonumber\\
	&\geq c_4\|(\Phi_y,\Psi_y,W_y)\|^2+c_4\v^{-1}\|(\phi_y,\psi_y,\zeta_y)\|^2-c_4\d\v^4(1+\v^2\tau)^{-\f32},
\end{align}
for some small positive constant $c_4$. Then, it follows from \eqref{N-2.105} and \eqref{N-2.106} that
\begin{equation}\label{N-2.107}
	\begin{cases}
		\|(\Phi_y,\Psi_y,W_y)\|^2+\|(\phi,\psi,\zeta)\|^2\leq C C_0\sqrt\d\v^2(1+\v^2\tau)^{-1+C_0\sqrt\d},\\
		\|(\phi_y,\psi_y,\zeta_y)\|^2\leq C C_0\sqrt\d\v^3(1+\v^2\tau)^{-1+C_0\sqrt\d}.
	\end{cases}
\end{equation}
From \eqref{N-2.102} and \eqref{N-2.107}, one has  the decay rate for $(\Phi,\Psi,W)$,
\begin{align}\label{N-2.108}
	\|(\Phi,\Psi,W)\|_{L^\infty} \leq C\|(\Phi,\Psi,W)\|^{\f12}\|(\Phi_y,\Psi_y,W_y)\|^{\f12}\leq C C_0\d^{\f14}\v^{\f12}(1+\v^2\tau)^{-\f14+\f12C_0\sqrt\d}.
\end{align}

However, from \eqref{N-2.95}, one observes that $\|(\phi_y,\psi_y,\zeta_y)\|^2$ should have better time-decay rate than the one in \eqref{N-2.107}. In fact, multiplying \eqref{N-2.95} by $(1+\v^2\tau)^{\f32}$, one gets that 
\begin{align}\label{N-2.109}
	[(1+\v^2\tau)^{\f32}E_{4}]_{\tau}&\leq C\d\v^3\|(\Phi_y,\Psi_y,W_y)\|^2+C\v(1+\v^2\tau)\|(\phi_y,\psi_y,\zeta_y)\|^2
	+C\d\v^6(1+\v^2\tau)^{-\f32}\nonumber\\
	&~~~+C(1+\v^2\tau)^{\f32}\|(\phi_y,\psi_y,\zeta_y)\|^{\f{10}{3}}
	+C(1+\v^2\tau)^{\f32}\|(\phi,\psi)\|\|(\phi_y,\psi_y,\zeta_y)\|^{3}.
\end{align}
By using \eqref{N-2.105} and \eqref{N-2.107}, one has that
\begin{align}\label{N-2.110}
&\int_{0}^{\tau}(1+\v^2s)^{\f32}\|(\phi_y,\psi_y,\zeta_y)\|^{\f{10}{3}}ds
+\int_{0}^{\tau}(1+\v^2s)^{\f32}\|(\phi,\psi)\|\|(\phi_y,\psi_y,\zeta_y)\|^{3}ds\nonumber\\
&\leq \v^2 \int_{0}^{\tau}(1+\v^2s)^{\f32}[(1+\v^2s)^{-\f23+\f23C_0\sqrt\d}
+(1+\v^2s)^{-1+C_0\sqrt\d}]\|(\phi_y,\psi_y,\zeta_y)\|^2ds\nonumber\\
&\leq \v^2 \int_{0}^{\tau}(1+\v^2s)\|(\phi_y,\psi_y,\zeta_y)\|^2ds
\leq C C_0\sqrt\d\v^4(1+\v^2\tau)^{C_0\sqrt\d},
\end{align}
for $C_0\sqrt\d\leq \f18$.
Thus, integrating \eqref{N-2.109} over $[0,\tau]$ and using \eqref{N-2.110}, \eqref{N-2.105} and \eqref{N-2.107}, one gets that 
\begin{equation*}\label{N-2.111}
	\|(\phi_y,\psi_y,\zeta_y)\|^2\leq CE_{4}\leq C C_0\sqrt\d\v^3(1+\v^2\tau)^{-\f32+C_0\sqrt\d}.
\end{equation*}
By using Sobolev inequality, we obtain
\begin{align}
	&\|(\phi,\psi,\zeta)\|_{L^\infty}\leq \|(\phi,\psi,\zeta)\|^{\f12}\|(\phi_y,\psi_y,\zeta_y)\|^{\f12}\nonumber\\
	&~~~~~~~~~~~~~~~~~~\leq C\d^\f14\v^{\f54}(1+\v^2\tau)^{-\f58+\f12C_0\sqrt\d}\leq C\d^\f14\v^{\f54}(1+\v^2\tau)^{-\f{9}{16}},\nonumber
\end{align}
for $C_0\sqrt\d\leq \f18$.
Therefore the {\it a priori } assumption \eqref{N-2.13} is closed and the proof of Proposition \ref{Npro1.1} is completed.  $\hfill\Box$

\section{Ill-prepared Data}

\subsection{Uniform Estimates}

For simplicity, throughout this subsection we omit  the superscript $\v$ of the variables and  
denote
\begin{equation}\label{2.10}
a(\v p)=e^{-\v p},～~b(\t)=e^{\t}.
\end{equation}
For later use, we define that
\begin{align}\label{2.3}
\begin{cases}
\mathcal{Q}(t):=\sup_{0\leq\tau\leq t}\Big\{\|(p,u)(\tau)\|_{\mathcal{H}^{s}}+\|(\v p, \v u)(\tau)\|_{\mathcal{H}^{s+1}}\\
~~~~~~~~~~~~~~~~~~~~~~~~~~~+\|(\t-\tilde{\t})(\tau)\|_{L^2}+\|((\v\partial_t) \t,\t_x)(\tau)\|_{\mathcal{H}^{s}}\Big\},\\[3mm]
\mathcal{S}(t):=\|(p_x,u_x)(t)\|_{\mathcal{H}^{s}}+\|(\v u_x,\t_x)(t)\|_{\mathcal{H}^{s+1}}.
\end{cases}
\end{align}
where we have used the notation $\|f(t)\|_{\mathcal{H}^{s}}:=\sum_{|\a|\leq s}\|\partial^\a f(t)\|_{L^2(\mathbb{R})}$ with $\partial^\a:=(\v\partial_t)^{\a_0}\partial_x^{\a_1}$ and $\a=(\a_0,\a_1)$.

We assume the following a priori assumption
\begin{equation}\label{2.14}
0<\f12\underline{a}\leq a(\v p)\leq 2\bar{a},~~0<\f12\underline{b}\leq b(\t)\leq 2\bar{b},~\mbox{on}~t\in[0,T].
\end{equation}

To prove Theorem \ref{thm1.1}, we  need only to  prove the following \textit{a priori} estimates. 
\begin{proposition}\label{pro2.1}
For any given integer $s\geq 4$ and  $\v\in(0,1]$, let $(p^\v,u^\v,\t^\v)$ be the classical solution to the Cauchy problem \eqref{2.1} and \eqref{2.1-1}. Under the a priori assumption \eqref{2.14}, it holds that 
\begin{equation*}
\mathcal{N}(t)\leq  C[1+\Lambda(\mathcal{Q}(0))]+C(t^{\f12}+\v)\Lambda(\mathcal{N}(t)),
\end{equation*}
where $\mathcal{N}(t)$ is defined in \eqref{2.2}, the constant $C>0$ may depend on $\underline{a}, \bar{a}, \underline{b}$ and $\bar{b}$, and $\Lambda(\cdot)$ is a finite order polynomial.
\end{proposition}

Firstly, we have

\begin{lemma}\label{lem2.1}
For $s\geq4$, it holds that
\begin{equation}\label{2.4}
\|\t-\tilde{\t}\|^2_{L^2}+\sum_{|\a|=1}^s\|\partial^\a \t\|^2_{L^2}+\sum_{|\a|=0}^s\int_0^t\|\partial^\a\t_x\|^2_{L^2}d\tau\leq C+\mathcal{Q}^2(0)+\int_0^t\Lambda(\mathcal{Q}(s))ds.
\end{equation}

\end{lemma}

\noindent\textbf{Proof}.  It follows from $\eqref{2.1}_3$ that
\begin{equation}\label{2.5}
(\t-\tilde{\t})_t+u(\t-\tilde{\t})_x+u_x-\k e^{-\v p}(e^{\t}(\t-\tilde{\t})_x)_x
=\tilde{\m}\v^2e^{-\v p}|u_x|^2-u\tilde{\t}_x+\k e^{-\v p}(e^\t \tilde{\t}_x)_x-\tilde{\t}_t.
\end{equation}
Multiplying \eqref{2.5} by $\t-\tilde{\t}$ and integrating over $\mathbb{R}$,  it holds that
\begin{align}\label{2.12}
&\f12\f{d}{dt}\int|\t-\tilde{\t}|^2dx+\k\int a(\v p)b(\t)|(\t-\tilde{\t})_x|^2dx+\int\k \v p_x e^{\t-\v p}(\t-\tilde{\t})(\t-\tilde{\t})_xdx\nonumber\\
&\leq C\|u_x\|_{L^\infty}\|\t-\tilde{\t}\|^2_{L^2}+ C\|u_x\|_{L^2}\|\t-\tilde{\t}\|_{L^2}+C\v^2\|\t-\tilde{\t}\|_{L^\infty}\|u_x\|^2_{L^2}\nonumber\\
&~~+\|\tilde{\t}_x\|_{L^\infty}(\|u\|_{L^2}+\|(\t-\tilde{\t})_x\|_{L^2})\|\t-\tilde{\t}\|_{L^2}+C\|\t-\tilde{\t}\|_{L^2}\nonumber\\
&\leq C+\Lambda(\mathcal{Q}(t)).
\end{align}
Integrating \eqref{2.12} with respect to time, one obtains that
\begin{align}\label{2.13}
\|(\t-\tilde{\t})(t)\|^2_{L^2}+\int_0^t\|(\t-\tilde{\t})_x(\tau)\|^2_{L^2}d\tau\leq C+\|(\t-\tilde{\t})(0)\|^2_{L^2}+\int_0^t\Lambda(\mathcal{Q}(\tau))d\tau.
\end{align}
Let $\t_\a=\partial^\a\t$, for $1\leq|\a|\leq s$. Applying $\partial^\a$ to $\eqref{2.1}_3$, one obtains that
\begin{align}\label{2.6}
&\partial_t\t_\a+u\partial_x\t_\a+\partial^\a u_x-\k e^{-\v p}(e^{\t}\partial_x\t_\a)_x\nonumber\\
&=-[\partial^\a,u]\t_x+\tilde{\m}\v^2\partial^\a(e^{-\v p^\v}|u^\v_x|^2)+\k\Big\{\partial^\a( e^{-\v p}(e^{\t}\t_x)_x)- e^{-\v p}(e^{\t}\partial_x\t_\a)_x\Big\}.
\end{align}
Multiplying \eqref{2.6} by  $\t_\a$ and integrating the resulting equation over $y$, one has that
\begin{align}\label{2.15}
&\f12\f{d}{dt}\int|\t_\a|^2dx
+\int a(\v p)b(\t)|\partial_x\t_\a|^2dx+\f{1}{2}\int u\partial_x(|\t_\a|^2)dx+\int \t_\a \partial^\a u_x dx\nonumber\\
&\leq C \Big|\int a(\v p) b(\t)\v p_x\t_\a\partial_x\t_\a dx\Big|+\|\t_\a\|_{L^2}\|[\partial^\a,u]\t_x\|_{L^2}+C\|\t_\a\|_{L^2}\|\partial^\a(e^{-\v p^\v}|\v u^\v_x|^2)\|_{L^2}\nonumber\\
&~~~+C\Big|\int \t_\a \Big[\partial^\a( e^{-\v p}(e^{\t}\t_x)_x)- e^{-\v p}(e^{\t}\partial_x\t_\a)_x\Big]dx\Big|.
\end{align}
Integrating by parts and using the Cauchy inequality, we obtain that 
\begin{align}\label{2.16}
&\Big|\int u\partial_x(|\t_\a|^2)dx\Big|+\Big|\int \t_\a \partial^\a u_x dx\Big|\leq \|u_x\|_{L^\infty}\|\t_\a\|^2_{L^2}+C\|\partial^\a u\|^2_{L^2}+C\|\partial_x\t_\a\|_{L^2}^2\nonumber\\
&\leq \|u\|_{H^2}\|\t_\a\|^2_{L^2}+C\|\partial^\a u\|^2_{L^2}+C\|\partial_x\t_\a\|_{L^2}^2\leq \Lambda(\mathcal{Q}(t)).
\end{align}
We shall estimate the RHS of \eqref{2.15}. Firstly, we have that
\begin{align}\label{2.17}
\Big|\int a(\v p) b(\t)\v p_x\t_\a\partial_x\t_\a dx\Big|\leq C\|\v p\|^2_{H^2}\|\t_\a\|^2_{L^2}+C\|\partial_x\t_\a\|^2_{L^2}\leq \Lambda(\mathcal{Q}(t)).
\end{align}
Notice that
\begin{align}\label{2.19}
[\partial^\a,u]\t_x=\sum_{1\leq|\b|\leq \a} C_{\a,\b}\partial^\b u\cdot \partial^{\a-\b}\t_x,
\end{align}
then one gets immediately  that
\begin{align}\label{2.20}
\|[\partial^\a,u]\t_x\|_{L^2}\leq C \Big(\sum_{|\b|\leq s}\|\partial^\b u\|_{L^2}\Big)\cdot\Big(\sum_{|\g|\leq s-1}\|\partial^\g \t_x\|_{L^2}\Big)\leq C\Lambda(\mathcal{Q}(t)).
\end{align}
A directly calculation shows that
\begin{align}\label{2.21}
\|\partial^\a(e^{-\v p^\v}|\v u^\v_x|^2)\|_{L^2}&\leq C\|\partial^\a(|\v u^\v_x|^2)\|_{L^2}+\Lambda(\|\v p\|_{\mathcal{H}^{s}})\cdot\|(\v u_x)^2\|_{\mathcal{H}^{s-1}}\nonumber\\
&\leq C\|\v u^\v_x\|^2_{\mathcal{H}^s}+\Lambda(\|\v p\|_{\mathcal{H}^{s}})\cdot\|(\v u_x)^2\|_{\mathcal{H}^{s-1}}\leq C\Lambda(\mathcal{Q}(t)).
\end{align}
Noting that
\begin{align}
&\partial^\a( e^{-\v p}(e^{\t}\t_x)_x)- e^{-\v p}(e^{\t}\partial_x\t_\a)_x\nonumber\\
&=\sum_{1\leq|\b|\leq\a}C_{\a,\b}\partial^\b(e^{-\v p})\cdot\partial^{\a-\b}(e^\t \t_x)_x
+e^{-\v p}\Big(\sum_{1\leq|\b|\leq\a}C_{\a,\b}\partial(e^{\t})\cdot\partial^{\a-\b}\t_x\Big)_x,\nonumber
\end{align}
which yields immediately that 
\begin{align}\label{2.22}
&\|\partial^\a( e^{-\v p}(e^{\t}\t_x)_x)- e^{-\v p}(e^{\t}\partial_x\t_\a)_x\|_{L^2}\nonumber\\
\leq &C\sum_{1\leq|\b|\leq\a} \Big\{\|\partial^\b(e^{-\v p})\cdot\partial^{\a-\b}(e^\t \t_{xx}+e^\t |\t_x|^2)\|_{L^2}+\|\partial^\b(e^\t \t_x)\cdot\partial^{\a-\b}\t_x\|\nonumber\\
&~~~~~~~~~~~~~~~~~~~~~~~~~~~~~~~~~~~~~~~+\|\partial^\b(e^\t)\cdot\partial^{\a-\b}\t_{xx}\|\Big\}\nonumber\\
\leq &C\Lambda(\|\v p\|_{\mathcal{H}^{s}})\cdot\Big(\|\t_{xx}\|_{\mathcal{H}^{s-1}}+\Lambda(\|\t_x\|_{\mathcal{H}^s})\Big)
+C\Lambda(\|\t_x\|_{\mathcal{H}^s})\cdot[1+\|\t_{xx}\|_{\mathcal{H}^{s-1}}]
\leq C\Lambda(\mathcal{Q}(t)).
\end{align}
Substituting \eqref{2.16}, \eqref{2.17}, \eqref{2.20}, \eqref{2.21} and \eqref{2.22} into \eqref{2.15}, one obtains that
\begin{align}\label{2.23}
&\f12\f{d}{dt}\int|\t_\a|^2dx
+\int a(\v p)b(\t)|\partial_x\t_\a|^2dx\leq  C\Lambda(\mathcal{Q}(t)).
\end{align}
Integrating \eqref{2.23} with respect to time, one gets that
\begin{align}\label{2.24}
\|\t_\a\|^2_{L^2}
+\int_0^t\|\partial_x\t_\a(\tau)\|^2d\tau\leq C\|\t_\a(0)\|^2_{L^2} +C\int_0^t\Lambda(\mathcal{Q}(\tau))d\tau.
\end{align}
Then, from  \eqref{2.24} and \eqref{2.13}, we obtain \eqref{2.4}. Thus the proof of Lemma \ref{lem2.1} is completed. $\hfill\Box$

\

\begin{lemma}\label{lem2.2}
For $s\geq4$, it holds that
\begin{align}\label{2.25}
\|(\v p,\v u)(t)\|^2_{\mathcal{H}^s}+\int_0^t\|\v u_x(\tau)\|^2_{\mathcal{H}^s}d\tau \leq C\mathcal{Q}^2(0)+C\int_0^t\Lambda(\mathcal{Q}(\tau))\cdot[1+\mathcal{S}(\tau)]d\tau.
\end{align}
\end{lemma}

\noindent\textbf{Proof}.  For simplicity, we denote
\begin{align}\label{2.26}
\breve{p}=\v p,~\breve{u}=\v u,~~\mbox{and}~~
(\breve{p}_\a,\breve{u}_\a)=\partial^\a(\breve{p},\breve{u}),~~0\leq|\a|\leq s.
\end{align}
Then the functions $(\breve{p},\breve{u})$ satisfy 
\begin{equation}\label{2.27}
\begin{cases}
\partial_t\breve{p}+u\breve{p}_x+(2u-a(\breve{p})b(\t)\t_x)_x
=a(\breve{p})|\breve{u}_x|^2+a(\breve{p})b(\t)\t_x\cdot\breve{p}_x,\\[2mm]
b(-\t)[\partial_t\breve{u}+u \breve{u}_x]+p_x=\tilde{\m}a(\breve{p})\breve{u}_{xx}.
\end{cases}
\end{equation}
Applying $\partial^\a$ to \eqref{2.27}, one gets that
\begin{align}\label{2.28}
\begin{cases}
\partial_t\breve{p}_\a+u\partial_x\breve{p}_{\a}=h_1+h_2+h_3+h_4,\\[2mm]
b(-\t)[\partial_t\breve{u}_\a+u\partial_x\breve{u}_{\a}]+\partial^\a p_x-\tilde{\m}a(\breve{p})\partial_{xx}\breve{u}_{\a}=h_5+h_6+h_7,
\end{cases}
\end{align}
where
\begin{align}\label{2.29}
\begin{cases}
& h_1=-[\partial^\a,u]\partial_x\breve{p},~~ h_2=-\partial^\a\partial_x\Big(2u-a(\breve{p})b(\t)\t_x\Big),\\
& h_3=\partial^\a\Big( a(\breve{p})|\breve{u}_x|^2\Big),~~ h_4=\partial^\a\Big(a(\breve{p})b(\t)\t_x\cdot\breve{p}_x\Big),\\
& h_5=-[\partial^\a, b(-\t)]\partial_t\breve{u},~~h_6=-[\partial^\a, b(-\t) u]\partial_x\breve{u},~~ h_7=\tilde{\m}[\partial^\a, a(\breve{p})]\breve{u}_{xx}.
\end{cases}
\end{align}
Multiplying $\eqref{2.28}_1$ by $\breve{p}_\a$, one has  that
\begin{align}\label{2.30}
\f12\f{d}{dt}\int|\breve{p}_\a|^2dx=\f12\int u_x|\breve{p}_\a|^2dx+\int(h_1+h_2+h_3+h_4)\breve{p}_\a dx.
\end{align}
It is straightforward  to imply that
\begin{equation}\label{2.31}
\Big|\int u_x|\breve{p}_\a|^2dx \Big|\leq \|u_x\|_{L^\infty}\|\breve{p}_\a\|^2_{L^2}\leq C\Lambda(\mathcal{Q}(t)).
\end{equation}
A directly calculation shows that
\begin{equation*}
\|(h_1,h_3,h_4)\|_{L^2}\leq C\Lambda(\mathcal{Q}(t)),
\end{equation*}
which yields immediately that
\begin{align}\label{2.32}
\Big|\int(h_1+h_3+h_4)\breve{p}_\a dx\Big|\leq C\Lambda(\mathcal{Q}(t)).
\end{align}
Noting that
\begin{align}\label{2.33}
\|h_2\|_{L^2}\leq C[\|\partial^\a u_x\|_{L^2}+\|\partial^\a\partial_{xx}\t\|_{L^2}+\Lambda(\mathcal{Q}(t))]
\leq C[\mathcal{S}(t)+\Lambda(\mathcal{Q}(t))],
\end{align}
we obtain that 
\begin{align}\label{2.34}
\Big|\int h_2\breve{p}_\a dx\Big|\leq C\|\breve{p}_\a\|_{L^2}\cdot[\mathcal{S}(t)+\Lambda(\mathcal{Q}(t))].
\end{align}
Substituting \eqref{2.31}, \eqref{2.32} and \eqref{2.34} into \eqref{2.30}, then integrating the resulting inequality with respect to time,  then we get that
\begin{align}\label{2.35}
\|\breve{p}_\a\|^2\leq \|\breve{p}_\a(0)\|^2+C\int_0^t\Lambda(\mathcal{Q}(\tau))\cdot[1+\mathcal{S}(\tau)]d\tau.
\end{align}

\

On the other hand, multiplying  $\eqref{2.28}_2$ by $\breve{u}_\a$ and integrating the resulting equation,  one obtains that
\begin{align}\label{2.36}
&\f12\f{d}{dt}\int b(-\t)|\breve{u}_\a|^2dx-\tilde{\m}\int a(\breve{p})\breve{u}_\a\partial_{xx}\breve{u}_\a dx\leq \int(h_5+h_6+h_7)\breve{u}_\a dx,
\end{align}
where we have used the facts 
\begin{align}\label{2.37}
\Big|\int \partial_t b(-\t)|\breve{u}_\a|^2dx\Big|+\Big|\int b(-\t)u(|\breve{u}_\a|^2)_xdx\Big|\leq C\Lambda(\mathcal{Q}(t)),
\end{align}
and
\begin{align}\label{2.38}
\Big|\int \partial^\a p_x\cdot \breve{u}_\a dx\Big|\leq \|\partial^\a p\|_{L^2}\cdot\|\partial_x\breve{u}_\a\|_{L^2} \leq C\Lambda(\mathcal{Q}(t)).
\end{align}
For the second term on the left hand side of \eqref{2.36}, it follows from integrating by parts that
\begin{align}\label{2.39}
-\tilde{\m}\int a(\breve{p})\breve{u}_\a\partial_{xx}\breve{u}_\a dx
&=\tilde{\m}\int a(\breve{p})|\partial_{x}\breve{u}_\a|^2 dx+\tilde{\m}\int \partial_xa(\breve{p})\breve{u}_\a\partial_{x}\breve{u}_\a dx\nonumber\\
&\geq \f34\tilde{\m}\int a(\breve{p})|\partial_{x}\breve{u}_\a|^2 dx-C\Lambda(\mathcal{Q}(t)).
\end{align}
Now we estimate the RHS of \eqref{2.36}. Noting 
\begin{align}
\|h_5\|_{L^2}=\|[\partial^\a,b(-\t)]\partial_t\breve{u}\|_{L^2}=\|[\partial^\a,b(-\t)](\v\partial_t)u\|_{L^2}\leq C\Lambda(\mathcal{Q}(t))\|u\|_{\mathcal{H}^s}\leq C\Lambda(\mathcal{Q}(t)),\nonumber
\end{align}
and
\begin{align}
\|h_6\|_{L^2}+\|h_7\|_{L^2}\leq C\Lambda(\mathcal{Q}(t))+C\Lambda(\mathcal{Q}(t))\|\v u_x\|_{\mathcal{H}^s}\leq C\Lambda(\mathcal{Q}(t)),\nonumber
\end{align}
which, together with Holder inequality, yield that 
\begin{align}\label{2.40}
\Big|\int(h_5+h_6+h_7)\breve{u}_\a dx\Big|\leq \Lambda(\mathcal{Q}(t)).
\end{align}
Substituting \eqref{2.40}-\eqref{2.37} into \eqref{2.36} and integrating the resulting inequality with respect to time, one obtains that
\begin{align}\label{2.41}
\|\v u(t)\|^2_{\mathcal{H}^s}+\int_0^t\|\v u(\tau)\|^2_{\mathcal{H}^s}d\tau
\leq C\|\v u(0)\|^2_{\mathcal{H}^s}+\int_0^t\Lambda(\mathcal{Q}(\tau))d\tau.
\end{align}
Combining \eqref{2.41} and \eqref{2.35}, one proves \eqref{2.25}.  Thus the proof of Lemma \ref{lem2.2} is completed. $\hfill\Box$

\

The following Lemma is devoted to the highest order derivative estimates $\|\partial^\a(\v p, \v u,\t)(t)\|^2_{L^2},$ $|\a|=s+1$.
\begin{lemma}\label{lem2.3}
For $s\geq 4$, it holds that 	
\begin{align}\label{2.44}
&\sum_{|\a|=s+1}\|\partial^\a(\v p, \v u,\t)(t)\|^2_{L^2}
+\sum_{|\a|=s+1} \int_0^t\|\partial^\a(\v u_x,\t_x)(\tau)\|^2_{L^2}d\tau\nonumber\\
&\leq C\mathcal{Q}^2(0)+C\int_0^t[1+\mathcal{S}(\tau)]\Lambda(\mathcal{Q}(\tau))d\tau.
\end{align}	
\end{lemma}

\noindent\textbf{Proof}. Similar to \cite{Alazard2,JJLX},  we set
\begin{align}\label{2.42}
(\check{p},\check{u})=(\v p-\t,\v u).
\end{align}
Then, from \eqref{2.1},  it is easy to check that  $(\check{p},\check{u},\t)$ satisfy 
\begin{align}\label{2.43}
\begin{cases}
\partial_t\check{p}+u \check{p}_x+\f1\v \check{u}_x=0,\\
b(-\t)[\partial_t\check{u}+u \check{u}_x]+\f1\v(\check{p}_x+\t_x)=\tilde{\m}a(\v p)\check{u}_{xx},\\
\partial_t\t+u \t_x+\f1\v \check{u}_x=a(\v p)(b(\t)\t_x)_x+\tilde{\m}\v^2 a(\v p)|u_x|^2.
\end{cases}
\end{align}
Let $\a$ be a multi-index with $|\a|=s+1$ and denote
\begin{align}\label{2.45}
(\check{p}_\a, \check{u}_\a, \t_\a)=\partial^\a(\check{p}, \check{u}, \t).
\end{align}
Applying $\partial^\a$ to \eqref{2.43}, one gets that
\begin{align}\label{2.46}
\begin{cases}
\partial_t\check{p}_\a+u \partial_x\check{p}_\a+\f1\v \partial_x\check{u}_\a=h_8,\\[2mm]
b(-\t)[\partial_t\check{u}_\a+u \partial_x\check{u}_\a]+\f1\v(\partial_x\check{p}_\a+\partial_x\t_\a)=\tilde{\m}a(\v p)\partial_{xx}\check{u}_{\a}+h_9,\\[2mm]
\partial_t\t_\a+u \partial_x\t_\a+\f1\v \partial_x\check{u}_\a=a(\v p)(b(\t)\partial_x\t_\a)_x+h_{10},
\end{cases}
\end{align}
where
\begin{align}
& h_8=-[\partial^\a,u]\check{p}_x,\nonumber\\
& h_9=-[\partial^\a,b(-\t)]\partial_t(\v u)-[\partial^\a,b(-\t)u]\partial_x(\v u)+[\partial^\a,\tilde{\m}a(\v p)]\partial_{xx}\check{u},\nonumber\\
& h_{10}=-[\partial^\a,u]\t_x+\tilde{\m}\partial^\a\Big(\v^2a(\v p)|u_x|^2\Big)+\Big\{\partial^\a(a(\v p)(b(\t)\t_x)_x)-a(\v p)(b(\t)\partial_x\t_\a)_x\Big\}.\nonumber
\end{align}
Considering $\int \big[\eqref{2.46}_1\times \check{p}_\a+\eqref{2.46}_2\times \check{u}_\a+\eqref{2.46}_1\times \t_\a \big] dx$ and integrating by parts the resulting equation, one can  obtain
\begin{align}\label{2.47}
&\f12\f{d}{dt}\int |\check{p}_\a|^2+b(-\t)|\check{u}_\a|^2+|\t_\a|^2dx
+\f34\int \tilde{\m}a(\v p)|\partial_{x}\check{u}_{\a}|^2+ a(\v p)b(\t)|\partial_x\t_\a|^2dx\nonumber\\
&\leq C\Lambda(\mathcal{Q}(t))+\int (h_8\check{p}_\a+h_9\check{u}_\a+h_{10}\t_\a)dx,
\end{align}
where we have used the following estimates
\begin{equation}\label{2.48}
\Big|\int u \check{p}_\a\partial_x\check{p}_\a dx\Big|+\Big|\int\partial_tb(-\t)\cdot|\check{u}_\a|^2dx\Big|
+\Big|\int b(\v p) u \check{u}_\a\partial_x\check{u}_\a dx\Big|+\Big|\int u \t_\a\partial_x\t_\a dx\Big|\leq \Lambda(\mathcal{Q}(t)),
\end{equation}
\begin{align}\label{2.50}
\int a(\v p)\partial_{xx}\check{u}_{\a}\check{u}_{\a}dx
&=-\int a(\v p)|\partial_{x}\check{u}_{\a}|^2dx-\int \partial_xa(\v p)\partial_{x}\check{u}_{\a}\check{u}_{\a}dx\nonumber\\
&\leq -\f34\int a(\v p)|\partial_{x}\check{u}_{\a}|^2dx+C\Lambda(\mathcal{Q}(t)),
\end{align}
\begin{align}\label{2.51}
\int a(\v p)(b(\t)\partial_{x}\t_{\a})_x\t_{\a}dx
&=-\int a(\v p)b(\t)|\partial_{x}\t_{\a}|^2dx-\int \partial_xa(\v p)\cdot b(\t)\t_{\a}\partial_{x}\t_{\a}dx\nonumber\\
&\leq -\f34\int a(\v p)b(\t)|\partial_{x}\check{u}_{\a}|^2dx+C\Lambda(\mathcal{Q}(t)),
\end{align}
and the fact that 
\begin{equation}\label{2.49}
	\int \left(\partial_x\check{u}_\a\check{p}_\a+\partial_x\check{p}_\a\check{u}_\a+\check{u}_\a\partial_x\t_\a
	+\partial_x\check{u}_\a\t_\a\right) dx=0.
\end{equation}

It remains to estimate the RHS of \eqref{2.47}. Firstly, it follows from $\eqref{2.1}_2$ that
\begin{align}\nonumber
\v\partial_tu=-u\partial_x(\v u)-b(\t)\partial_xp+\tilde{\m} a(\v p)b(\t)\partial_{xx}(\v u),
\end{align}
which yields immediately that
\begin{equation}\label{2.52}
\|\v\partial_tu\|_{\mathcal{H}^s}\leq C\Lambda(\mathcal{Q}(t))\Big\{1+\|\partial_x(\v u)\|_{\mathcal{H}^s}+\|\partial_xp\|_{\mathcal{H}^s}+\|\partial_{xx}(\v u)\|_{\mathcal{H}^s}\Big\}\leq C\Lambda(\mathcal{Q}(t))\Big\{1+\mathcal{S}(t)\Big\}.
\end{equation}
Using \eqref{2.52}, one gets that
\begin{align}
\|h_8\|_{L^2}&=\|[\partial^\a,u](\v p_x-\t_x)\|_{L^2}\leq C\|\v p_x-\t_x\|_{\mathcal{H}^{s}}\|u\|_{\mathcal{H}^{s+1}}\nonumber\\
&\leq \Lambda(\mathcal{Q}(t))\Big[\|\v\partial_tu\|_{\mathcal{H}^s}+\|u_x\|_{\mathcal{H}^s}+\|u\|_{\mathcal{H}^s}\Big]\nonumber\\
&\leq  C\Lambda(\mathcal{Q}(t))\Big\{1+\mathcal{S}(t)\Big\},\nonumber
\end{align}
which, together with Holder inequality, implies that 
\begin{align}\label{2.53}
\Big|\int h_8 \check{p}_\a dx\Big|\leq  C\Lambda(\mathcal{Q}(t))\Big\{1+\mathcal{S}(t)\Big\}.
\end{align}
Next, we shall estimate the terms involving $h_9$. It is straightforward to imply that
\begin{align}\label{2.54}
&\|[\partial^\a,b(-\t)u]\partial_x(\v u)\|_{L^2}+\|[\partial^\a,\tilde{\m}a(\v p)]\partial_{xx}\check{u}\|_{L^2}\nonumber\\
&\leq C\Lambda(\mathcal{Q}(t))\Big[1+\|\v u\|_{\mathcal{H}^{s+1}}+\|\v u_{xx}\|_{\mathcal{H}^{s}}\Big]
\leq C\Lambda(\mathcal{Q}(t))\Big[1+\mathcal{S}(t)\Big],
\end{align}
and
\begin{align}\label{2.55}
&\|[\partial^\a,b(-\t)]\partial_t(\v u)\|_{L^2}\leq C\Lambda(\mathcal{Q}(t))[1+\|\v\partial_tu\|_{\mathcal{H}^{s}}]
\leq C\Lambda(\mathcal{Q}(t))\Big[1+\mathcal{S}(t)\Big],
\end{align}
where we have used \eqref{2.52} in the last inequality of  \eqref{2.55}. Then, it follows from \eqref{2.54} and \eqref{2.55} that
\begin{align}\label{2.57}
\Big|\int h_9 \check{u}_\a dx \Big|\leq C\Lambda(\mathcal{Q}(t))\Big[1+\mathcal{S}(t)\Big].
\end{align}
Finally, we estimate the terms involving $h_{10}$. Similarly, we have that
\begin{align}\nonumber
\|h_{10}\|_{L^2}\leq C\Lambda(\mathcal{Q}(t))\Big[1+\mathcal{S}(t)\Big],
\end{align}
which yields  that
\begin{align}\label{2.58}
\Big|\int h_{10} \t_\a dx \Big|\leq C\Lambda(\mathcal{Q}(t))\Big[1+\mathcal{S}(t)\Big].
\end{align}
Substituting \eqref{2.58}, \eqref{2.57} and \eqref{2.53} into \eqref{2.47} and integrating the resulting inequality with respect to time, one gets \eqref{2.44}. Thus the proof of Lemma \ref{lem2.3} is completed.   $\hfill\Box$

\

To establish the estimates for $\|(p,u)\|_{\mathcal{H}^s}$, we first control the term $(\v\partial_t)^s(p,u)$ which plays key role in the estimates for  $\|(p,u)\|_{\mathcal{H}^s}$. We start with a $L^2$-estimate for the following linearized equations around a given state $(\underline{p},\underline{u}, \underline{\theta})$
\begin{align}\label{2.59}
\begin{cases}
p_t^l+\underline{u}\cdot p_x^l+\f{1}{\v}(2u^l-\k a(\v \underline{p})b(\underline{\t})\t^l_x)_x
=\tilde\m\v a(\v \underline{p})\underline{u}_{x}u^l_x+\k a(\v \underline{p})b(\underline{\t})\underline{p}_x\t^l_x +f_1,\\[2mm]
b(-\underline{\t})[u^l_t+\underline{u}\cdot u^l_x]+\f1\v p^l_x=\tilde{\mu}a(\v \underline{p}) u^l_{xx}+f_2,\\[2mm]
\t^l_t+\underline{u}\t^l_x+u^l_x=\k a(\v \underline{p})(b(\underline{\t})\t^l_x)_x+\tilde{\m}\v^2a(\v \underline{p})\underline{u}_x u^l_x+f_3,\\[1mm]
\end{cases}
\end{align}
where we $f_i,~i=1,2,3$ are source terms.

\begin{lemma}\label{lem2.4}
Let $(p^l,u^l,\t^l)$  be the solution of \eqref{2.59} and assume that
\begin{equation*}\label{2.14-2}
0<\f{1}{2}\underline{a}\leq a(\v\underline{p})\leq 2 \bar{a},~~\mbox{and}~~0<\f{1}{2}\underline{b}\leq b(\underline{\t})\leq 2\bar{b},~\mbox{on}~t\in[0,T],
\end{equation*}
then it holds that, for $0<t\le T$,
\begin{align}\label{2.61}
&\|(p^l, u^l)(t)\|^2_{L^2}+\int_0^t\|u^l_x(\tau)\|^2_{L^2}d\tau\nonumber\\
&\leq C\|(p^l, u^l)(0)\|^2_{L^2}+C\sup_{0\leq\tau\leq t}\|\t^l_x(\tau)\|^2_{L^2}+C\Lambda(R_0)\bigg\{\int_0^t\|(\v u^l_x, \t^l_{xx})\|^2_{L^2}d\tau\\
& ~~+\int_0^t\|(p^l, u^l, \t^l_{x})\|^2_{L^2}d\tau+\int_0^t\|f_3\|^2_{L^2}d\tau
+\big(\int_0^t\|(p^l, u^l , \t^l_{x})\|^2_{L^2}d\tau\big)^{\f12}\cdot\big(\int_0^t\|f_1,f_2\|^2_{L^2}d\tau\big)^{\f12}\bigg\},\nonumber
\end{align}
where $R_0\doteq \sup_{\tau\in[0,T]}\{\|\partial_t\underline{\t}(\tau)\|_{L^\infty}+\|(\underline{p},\underline{u},\underline{\t})(\tau)\|_{W^{1,\infty}}\}$.
\end{lemma}

\noindent\textbf{Proof}.  Following \cite{JJLX}, we define
\begin{align}\label{2.62}
\mathfrak{u} := 2u^l-\k a(\v\underline{p})b(\underline{\t})\t^l_x.
\end{align}
Then,  it is straightforward to check that $p^l$ solves 
\begin{align}\label{2.63}
p^l_t+\underline{u}\cdot p^l_x+\f{1}{\v}\mathfrak{u}_x
=\tilde\m\v a(\v \underline{p})\underline{u}_{x}\mathfrak{u}_x+\f{\tilde\m}2\v a(\v \underline{p})\underline{u}_{x}(\k a(\v \underline{p})b(\underline{\t})\t^l_x)_x+\k a(\v \underline{p})b(\underline{\t})\underline{p}_x\cdot\t^l_x +f_1.
\end{align}
Consider $\f12 a(\v \underline{p})\cdot\partial_x\eqref{2.59}_3$, one can get that 
\begin{align}
&\f12 b(-\underline{\t})\Big[\partial_t(a(\v\underline{p})b(\underline{\t})\t^l_x)+\underline{u}\cdot(a(\v\underline{p})b(\underline{\t})\t^l_x)_x\Big]\nonumber\\
&=\f12 b(-\underline{\t})\Big[\partial_t(a(\v\underline{p})b(\underline{\t}))\t^l_x+\underline{u}\cdot(a(\v\underline{p})b(\underline{\t}))_x\t^l_x\Big]-\f{\k}{4} a(\v \underline{p})[a(\v\underline{p})b(\underline{\t})\t^l_x]_{xx}\nonumber\\
&~~~+\f{\tilde{\m}}4 \v^2a(\v \underline{p})[a(\v \underline{p})\underline{u}_x\mathfrak{u}_x]_x
+\f{\tilde{\m}}4 \v^2a(\v\underline{p})[a(\v\underline{p})\underline{u}_x\cdot(a(\v\underline{p})b(\underline{\t})\t^l_x)_x]_x
-\f12 a(\v \underline{p})\underline{u}_x\t^l_x\nonumber\\
&~~~+\f12 a(\v \underline{p})[\k a(\v \underline{p})(b(\underline{\t})\t^l_x)_x]_x-\f14 a(\v \underline{p})\mathfrak{u}_{xx}+\f12 a(\v \underline{p})\cdot \partial_xf_3,\nonumber
\end{align}
which, together with $\eqref{2.59}_2$, yields that 
\begin{align}\label{2.65}
&\f12 b(-\underline{\t})\Big[\partial_t\mathfrak{u}+\underline{u}\cdot\mathfrak{u}_x\Big]+\f1\v p^l_x-\f14 a(\v \underline{p})\mathfrak{u}_{xx}-\f12\tilde{\m}a(\v \underline{p})\mathfrak{u}_{xx}\nonumber\\
&=-\f12 b(-\underline{\t})\Big[\partial_t(a(\v\underline{p})b(\underline{\t}))\t^l_x+\underline{u}\cdot(a(\v\underline{p})b(\underline{\t}))_x\t^l_x\Big]+\f{\k}{4} a(\v \underline{p})[a(\v\underline{p})b(\underline{\t})\t^l_x]_{xx}\nonumber\\
&~~~-\f{\tilde{\m}}4 \v^2a(\v \underline{p})[a(\v \underline{p})\underline{u}_x\mathfrak{u}_x]_x
-\f{\tilde{\m}}4 \v^2a(\v\underline{p})[a(\v\underline{p})\underline{u}_x\cdot(a(\v\underline{p})b(\underline{\t})\t^l_x)_x]_x
+\f12 a(\v \underline{p})\underline{u}_x\t^l_x\nonumber\\
&~~~-\f12 a(\v \underline{p})[\k a(\v \underline{p})(a(\underline{\t})\t^l_x)_x]_x+\f{\tilde{\m}}2 a(\v \underline{p})\cdot[\k a(\v\underline{p})b(\underline{\t})\t^l_x]_{xx}+f_2-\f12 a(\v \underline{p})\cdot \partial_xf_3\nonumber\\
&=:g+f_2-\f12 a(\v \underline{p})\cdot \partial_xf_3.
\end{align}
Multiplying \eqref{2.63}  by $p^l$,   \eqref{2.65} by $\mathfrak{u}$, adding them together and integrating the resulting equation, we get 
\begin{align}\label{2.66}
&\f{d}{dt}\left(\int\f12|p^l|^2+\f14b(-\underline{\t})|\mathfrak{u}|^2dx  \right)
+\int (\f14+\f12\tilde{\m})a(\v \underline{p})|\mathfrak{u}_x|^2dx\nonumber\\
& = \f12\int \underline{u}_x|p^l|^2dx+ \f14\int [\partial_tb(-\underline{\t})+\partial_x(b(-\underline{\t})\underline{u})]|\mathfrak{u}|^2dx+\int (\f14+\f12\tilde{\m})\v a(\v \underline{p}) \underline{p}_x \mathfrak{u}_x \mathfrak{u} dx\nonumber\\
&~~~~+\int \tilde\m \v a(\v \underline{p})\underline{u}_{x}\mathfrak{u}_x p^l dx+\int\f{\tilde\m}2\v a(\v \underline{p})\underline{u}_{x}(\k a(\v \underline{p})b(\underline{\t})\t^l_x)_x p^l dx+\k \int a(\v \underline{p})b(\underline{\t})\underline{p}_x\cdot\t^l_x p^l dx\nonumber\\
&~~~~~+\int g \mathfrak{u}dx+\int (f_1 p^l+f_2\mathfrak{u})dx-\f12 \int a(\v \underline{p}) \partial_xf_3 \cdot\mathfrak{u}dx\nonumber\\
&\leq \f18\int (\f14+\f12\tilde{\m})a(\v \underline{p})|\mathfrak{u}_x|^2dx+C\Lambda(R_0)\|(p^l, u^l, \t^l_x,\t^l_{xx})\|^2+ \|f_3\|^2_{L^2}\nonumber\\
&~~~~~+\int (f_1 p^l+f_2 \mathfrak{u}) dx+\int g \mathfrak{u}dx.
\end{align}
It is easy to note that 
\begin{align}\label{2.67}
\Big|\int g \mathfrak{u}dx\Big|\leq \f18\int (\f14+\f12\tilde{\m})a(\v \underline{p})|\mathfrak{u}_x|^2dx+C\Lambda(R_0)\|(p^l,u^l,\t_x^l,\t_{xx}^l)\|_{L^2}^2.
\end{align}
Substituting \eqref{2.67} into \eqref{2.66} and integrating the resulting inequality with respect to time,  one obtains 
\eqref{2.61}.  Therefore we complete the proof of Lemma \ref{lem2.4}.   
$\hfill\Box$

\

Next we shall use Lemma \ref{lem2.4} to estimate  $\|(\v\partial_t)^k(p,u)\|_{L^2}$ for $1\leq k\leq s$.
\begin{lemma}\label{lem2.5}
For $s\geq4$ and $0\leq k\leq s$, it holds that 
\begin{align}\label{2.68}
&\|(\v\partial_t)^k(p,u)(t)\|^2_{L^2}+\int_0^t\|(\v\partial_t)^ku_x(\tau)\|^2_{L^2}d\tau\nonumber\\
&\leq C\|(\v\partial_t)^k(p,u)(0)\|^2_{L^2}+C\sup_{0\leq \tau\leq t}\|\t_x(\tau)\|^2_{\mathcal{H}^k}+Ct^{\f12}\Lambda(\mathcal{N}(t))+C\Lambda(R)\int_0^t\|\t_{xx}(\tau)\|^2_{\mathcal{H}^k}d\tau,
\end{align}
where $R(t)=\sup_{0\leq \tau\leq t}\Big\{\|\partial_t\t(\tau)\|_{L^\infty}+\|(p,u,\t)(\tau)\|_{W^{1,\infty}}\Big\}$.
\end{lemma}

\noindent\textbf{Proof}. For simplicity, we set 
\begin{align}\label{2.70}
(p_k,u_k,\t_k)=(\v\partial_t)^k(p,u,\t),~~0\leq k\leq s.
\end{align}
Applying $(\v\partial_t)^k$ to \eqref{2.1}, one has that
\begin{align}\label{2.69}
\begin{cases}
p_{kt}+u\cdot p_{kx}+\f{1}{\v}(2u_k-\k a(\v p)b(\t)\t_{kx})_x
=\tilde\m\v a(\v p)u_x u_{kx}+\k a(\v p)b(\t)p_x\cdot\t_{kx}+f_{k1},\\[2mm]
b(-\t)[u_{kt}+u\cdot u_{kx}]+\f1\v p_{kx}=\tilde{\mu}a(\v p)u_{kxx}+f_{k2},\\[2mm]
\t_{kt}+u\t_{kx}+u_{kx}=\k a(\v p)(b(\t)\t_{kx})_x+\tilde{\m}\v^2a(\v p)u_xu_{kx}+f_{k3},
\end{cases}
\end{align}
where 
\begin{align}\label{2.71}
\begin{cases}
 f_{k1}=-[(\v\partial_t)^k,u]p_x+\k\f1\v\Big([(\v\partial_t)^k,a(\v p)b(\t)]\t_{x}\Big)_x+\v[(\v\partial_t)^k,\tilde\m\v a(\v p)u_x]u_x\\
~~~~~~~~~~~~~~~~~~~~~~~+[(\v\partial_t)^k,\k a(\v p)b(\t)p_x]\t_x,\\[2mm]
 f_{k2}=-[(\v\partial_t)^k,b(-\t)]u_t-[(\v\partial_t)^k,b(-\t)u]u_x+\tilde{\m}[(\v\partial_t)^k,a(\v p)]u_{xx},\\[2mm]
 f_{k3}=-[(\v\partial_t)^k,u]\t_x+\tilde{\m}\v[(\v\partial_t)^k, a(\v p)u_x]u_x
+\k[(\v\partial_t)^k,a(\v p)](b(\t)\t_x)_x\\
~~~~~~~~~~~~~~~~~~~~~~~+\k a(\v p)([(\v\partial_t)^k,b(\t)]\t_x)_x.
\end{cases}
\end{align}
Applying  \ref{lem2.4} to equations \eqref{2.69}, we obtains  
\begin{align}\label{2.72}
&\|(\v\partial_t)^k(p,u)(t)\|^2_{L^2}+\int_0^t\|(\v\partial_t)^k u_x\|^2_{L^2}d\tau\nonumber\\
&\leq C\|(\v\partial_t)^k(p,u)(0)\|^2_{L^2}+C\sup_{0\leq \tau\leq t}\|\t_x(\tau)\|^2_{\mathcal{H}^k}
+C\Lambda(R)\bigg\{\int_0^t\|(\v\partial_t)^k\t_xx(\tau)\|^2_{L^2}d\tau\nonumber\\
&~~~~+\int^t_0\|(p_k, u_k, (\t_k)_x)\|^2_{L^2} d\tau+\int_0^t\|f_{k3}\|^2_{L^2}d\tau
+t^{\f12}\Lambda(\mathcal{Q}(t))\big(\int_0^t\|(f_{k1},f_{k2})\|^2_{L^2}d\tau\big)^{\f12}\bigg\}.
\end{align}

It remains to estimate the terms involving  $(f_{k1},f_{k2},f_{k3})$. We need only to estimate the case for $k\geq1$ since $(f_{01},f_{02},f_{03})=(0,0,0)$. 
For the singular term of $f_{k1}$, we first notice that 
\begin{align}\label{2.73}
\Big([(\v\partial_t)^k,a(\v p)b(\t)]\t_{x}\Big)_x=[(\v\partial_t)^k,a(\v p)b(\t)]\t_{xx}+[(\v\partial_t)^k,\partial_x(a(\v p)b(\t))]\t_{x},
\end{align}
then a careful calculation gives that 
\begin{align}\label{2.74}
&\f1\v\|[(\v\partial_t)^k,a(\v p)b(\t)]\t_{xx}\|_{L^2}=\|[(\v\partial_t)^{k-1}\partial_t,a(\v p)b(\t)]\t_{xx}\|_{L^2}\nonumber\\
&\leq \|[(\v\partial_t)^{k-1},a(\v p)b(\t)]\partial_t\t_{xx}\|_{L^2}+ \|[(\v\partial_t)^{k-1},\partial_t(a(\v p)b(\t))]\t_{xx}\|_{L^2}\nonumber\\
&\leq C\Lambda(\mathcal{Q}(t))[1+\|\partial_t\t_x\|_{\mathcal{H}^{s-1}}+\|\partial_t\t\|_{\mathcal{H}^{s-1}}],
\end{align}
and
\begin{align}\label{2.75}
&\f1\v\|[(\v\partial_t)^k,\partial_x(a(\v p)b(\t))]\t_{x}\|_{L^2}=\|[(\v\partial_t)^{k-1}\partial_t,\partial_x(a(\v p)b(\t))]\t_{x}\|_{L^2}\nonumber\\
&\leq \|[(\v\partial_t)^{k-1}\partial_t,\partial_{xt}(a(\v p)b(\t))]\t_{x}\|_{L^2}
+ \|[(\v\partial_t)^{k-1},\partial_x(a(\v p)b(\t))]\t_{xt}\|_{L^2}\nonumber\\
&\leq C\Lambda(\mathcal{Q}(t))[1+\|\partial_t\t_x\|_{\mathcal{H}^{s-1}}+\|p_x\|_{\mathcal{H}^{s}}].
\end{align}
For the remain terms of $f_{k1}$, it is straightforward to obtain
\begin{align}\label{2.76}
&\|[(\v\partial_t)^k,u]p_x\|_{L^2}+\|\v[(\v\partial_t)^k,\tilde\m\v a(\v p)u_x]u_x\|_{L^2}+\|[(\v\partial_t)^k,\k a(\v p)b(\t)p_x]\t_x\|_{L^2}\nonumber\\
&\leq C\Lambda(\mathcal{Q}(t))[1+\|p_x\|_{\mathcal{H}^{s}}].
\end{align}
Combining \eqref{2.73}-\eqref{2.76}, one gets that 
\begin{align}\label{2.77}
\|f_{k1}\|_{L^2}\leq C\Lambda(\mathcal{Q}(t))[1+\|\partial_t\t_x\|_{\mathcal{H}^{s-1}}+\|p_x\|_{\mathcal{H}^{s}}].
\end{align}

For the estimate of  $\|f_{k2}\|_{L^2}$,  it is noted that 
\begin{align}\nonumber
[(\v\partial_t)^k,b(-\t)]u_t=\sum_{i=1}^k C_{k,i}(\v\partial_t)^ib(-\t)\cdot (\v \partial_t)^{k-i}u_t
=\sum_{i=1}^k C_{k,i}(\v\partial_t)^{i-1}\partial_tb(-\t)\cdot (\v \partial_t)^{k-i+1}u,
\end{align}
which yields immediately that 
\begin{align}\label{2.78}
\|[(\v\partial_t)^k,b(-\t)]u_t\|_{L^2}\leq  C\Lambda(\mathcal{Q}(t))[1+\|\partial_t\t\|_{\mathcal{H}^{k-1}}+\Lambda(\|\t_t\|_{\mathcal{H}^{k-2}})], 
\end{align}
in which  $\|\t_t\|_{\mathcal{H}^{k-2}}$ does not appear when $k-2<0$. For the second and third terms of $f_{k2}$,  it is straightforward to obtain that 
\begin{align}\label{2.79}
\|[(\v\partial_t)^k,b(-\t)u]u_x\|_{L^2}+\|\tilde{\m}[(\v\partial_t)^k,a(\v p)]u_{xx}\|_{L^2}\leq C\Lambda(\mathcal{Q}(t)).
\end{align}
Then it follows from  \eqref{2.78} and \eqref{2.79} that 
\begin{align}\label{2.80}
\|f_{k2}\|_{L^2}\leq C\Lambda(\mathcal{Q}(t))[1+\|\partial_t\t\|_{\mathcal{H}^{k-1}}+\Lambda(\|\t_t\|_{\mathcal{H}^{k-2}})].
\end{align}
Finally, a direct calculation yields that 
\begin{align}\label{2.81}
\|f_{k3}\|_{L^2}\leq C\Lambda(\mathcal{Q}(t)).
\end{align}

Note that both  \eqref{2.77} and \eqref{2.80} contain some norms of  $\partial_t\t$, which are not included in $\mathcal{Q}(t)$. To bound $\partial_t\t$, we use  $\eqref{2.1}_3$  to obtain 
\begin{align}\label{2.82}
\|\t_{t}\|_{\mathcal{H}^{s-1}}\leq C\Lambda(\mathcal{Q}(t)),~~\mbox{and}~~\|\t_{tx}\|_{\mathcal{H}^{s-1}}\leq C\Lambda(\mathcal{Q}(t))\Big(1+\|u_x\|_{\mathcal{H}^{s}}+\|\t_{xx}\|_{\mathcal{H}^{s}}\Big).
\end{align}
Then, combining \eqref{2.82}, \eqref{2.80} and \eqref{2.77}, one gets that 
\begin{align}\label{2.83}
\|(f_{k1},f_{k2})\|_{L^2}\leq  C\Lambda(\mathcal{Q}(t))[1+\|(p_x,u_x,\t_{xx})\|_{\mathcal{H}^{s}}].
\end{align}
Substituting \eqref{2.81} and \eqref{2.83} into \eqref{2.72}, one proves  \eqref{2.68}. Therefore, the proof of Lemma \ref{lem2.5} is completed.
$\hfill\Box$

\

From Lemma \ref{lem2.1}, Lemma \ref{lem2.3} and Lemma \ref{lem2.5}, we have the following corollary:
\begin{corollary}\label{cor2.1}
For $s\geq4$ and $0\leq k\leq s-1$, it holds that
\begin{align}\label{2.84}
&\|(\v\partial_t)^k(p,u)(t)\|^2_{L^2}+\int_0^t\|(\v\partial_t)^ku_x(\tau)\|^2_{L^2}d\tau\leq C\Big(1+\mathcal{Q}^2(0)\Big)+Ct^{\f12}\Lambda(\mathcal{N}(t)),
\end{align}
and
\begin{equation}\label{2.85}
\|(\v\partial_t)^s(p,u)(t)\|^2_{L^2}+\int_0^t\|(\v\partial_t)^su_x(\tau)\|^2_{L^2}d\tau\leq C\Big(1+\Lambda(\mathcal{Q}(0))\Big)+Ct^{\f12}\Lambda(\mathcal{N}(t))+C\Lambda(R(t)).
\end{equation}
where $R(t)$ is defined in Lemma \ref{lem2.5}.
\end{corollary}

\

We now use Corollary \ref{cor2.1} to estimate $\|(p,u)\|_{\mathcal{H}^{s}}$.  It follows from \eqref{2.1} that 
\begin{align}\label{2.86}
\begin{cases}
2u_x=-\v\partial_t p- \v u  p_x+(\k a(\v p)b(\t)\t_{x})_x
+\tilde\m a(\v p)|\v u_{x}|^2+\k a(\v p)b(\t)\v p_x\cdot\t_{x},\\[2mm]
p_{x}=-b(-\t)\v\partial_tu-\v b(-\t)u u_{x}+\tilde{\mu}a(\v p)\v u_{xx}.
\end{cases}
\end{align}

\begin{lemma}\label{lem2.7}
It holds, for $s\geq4$,  that 
\begin{align}\label{2.96}
\|(p,u)\|_{\mathcal{H}^{s}}\leq C\Big(1+\Lambda(\mathcal{Q}(0))\Big)+C(t^{\f12}+\v)\Lambda(\mathcal{N}(t)).
\end{align}
\end{lemma}
\noindent\textbf{Proof}. It follows from $\eqref{2.86}_1$, Lemma \ref{lem2.1}-Lemma \ref{lem2.3} and  Corollary \ref{cor2.1}  that 
\begin{align}\label{2.87}
\|u_x\|_{L^2}&\leq C\Big\{\|(\v\partial_t) p\|_{L^2}+\v\|u\|_{H^1}\|p_x\|_{L^2}+\|\t_{xx}\|_{L^2}+\Lambda(\|(\v p,\v u_x)\|_{\mathcal{H}^s}+\|\t_x\|_{\mathcal{H}^{s-1}})\Big\}\nonumber\\
&\leq C\Big(1+\Lambda(\mathcal{Q}(0))\Big)+C(t^{\f12}+\v)\Lambda(\mathcal{N}(t))
\end{align}
Similarly, one has,  for $0\leq k\leq s-2$, that 
\begin{align}\label{2.88}
\|(\v\partial_t)^ku_x\|_{L^2}&\leq C\Big\{\|(\v\partial_t)^{k+1} p\|_{L^2}+\v\Lambda(\mathcal{Q})+\Lambda(\|(\v p,\v u)\|_{\mathcal{H}^s}+\|(\v \partial_t\t,\t_x,\t_{xx})\|_{\mathcal{H}^{s-1}})\Big\}\nonumber\\
&\leq C\Big(1+\Lambda(\mathcal{Q}(0))\Big)+C(t^{\f12}+\v)\Lambda(\mathcal{N}(t)),
\end{align}
and 
\begin{align}\label{2.89}
\|(\v\partial_t)^{s-1}u_x\|_{L^2}&\leq C\Big\{\|(\v\partial_t)^{s} p\|_{L^2}+\v\Lambda(\mathcal{Q})+\Lambda(\|(\v p,\v u)\|_{\mathcal{H}^s}+\|(\v \partial_t\t,\t_x,\t_{xx})\|_{\mathcal{H}^{s-1}})\Big\}\nonumber\\
&\leq C\Big(1+\Lambda(\mathcal{Q}(0))\Big)+C(t^{\f12}+\v)\Lambda(\mathcal{N}(t))+C\Lambda(R(t)).
\end{align}
For $p_x$, it follows from $\eqref{2.86}_2$, Lemma \ref{lem2.1}-Lemma \ref{lem2.3} and  Corollary \ref{cor2.1}, for $0\leq k\leq s-2$,   that 
\begin{align}\label{2.90}
\|(\v\partial_t)^kp_x\|_{L^2}&\leq C\Big\{\|(\v\partial_t)^{k+1} u\|_{L^2}+\v\Lambda(\mathcal{Q})+\Lambda(\|(\v p,\v u)\|_{\mathcal{H}^s}+\|(\v \partial_t\t,\v u_{xx})\|_{\mathcal{H}^{s-1}})\Big\}\nonumber\\
&\leq C\Big(1+\Lambda(\mathcal{Q}(0))\Big)+C(t^{\f12}+\v)\Lambda(\mathcal{N}(t)),
\end{align}
and 
\begin{align}\label{2.91}
\|(\v\partial_t)^{s-1}p_x\|_{L^2}&\leq C\Big\{\|(\v\partial_t)^{s} u\|_{L^2}+\v\Lambda(\mathcal{Q})+\Lambda(\|(\v p,\v u)\|_{\mathcal{H}^s}+\|(\v \partial_t\t,\v u_{xx})\|_{\mathcal{H}^{s-1}})\Big\}\nonumber\\
&\leq C\Big(1+\Lambda(\mathcal{Q}(0))\Big)+C(t^{\f12}+\v)\Lambda(\mathcal{N}(t))+C\Lambda(R(t)).
\end{align}
Using the system \eqref{2.86} and \eqref{2.88}, \eqref{2.90}, one obtains, for $0\leq k\leq s-3$ that 
\begin{align}\label{2.92}
\|(\v\partial_t)^{k}(p_{xx},u_{xx})\|_{L^2}\leq C\Big(1+\Lambda(\mathcal{Q}(0))\Big)+C(t^{\f12}+\v)\Lambda(\mathcal{N}(t)).
\end{align}
Similarly,   one can get that 
\begin{align}\label{2.93}
\|(p,u)\|_{\mathcal{H}^{s-1}}\leq C\Big(1+\Lambda(\mathcal{Q}(0))\Big)+C(t^{\f12}+\v)\Lambda(\mathcal{N}(t)).
\end{align}
Then it follows from $\eqref{2.1}_3$, \eqref{2.4} and  \eqref{2.93} that
\begin{align}\label{2.94}
R(t)\leq C\Big(1+\Lambda(\mathcal{Q}(0))\Big)+C(t^{\f12}+\v)\Lambda(\mathcal{N}(t)).
\end{align}
Combining \eqref{2.85} \eqref{2.86}, \eqref{2.89}, \eqref{2.91}, \eqref{2.93}, \eqref{2.94}, and using the same argument as in \eqref{2.93}, one can get  that  
\begin{align}\label{2.95}
\|(p,u)\|_{\mathcal{H}^{s}}\leq C\Big(1+\Lambda(\mathcal{Q}(0))\Big)+C(t^{\f12}+\v)\Lambda(\mathcal{N}(t)).
\end{align}
Therefore the proof of Lemma \ref{lem2.7} is completed.  
$\hfill\Box$

\

\begin{lemma}\label{lem2.8}
It holds, for $s\geq4$,  that 
\begin{align}\label{2.96-1}
\int_0^t\|(p_x,u_x)(\tau)\|^2_{\mathcal{H}^s}d\tau\leq C\Big(1+\Lambda(\mathcal{Q}(0))\Big)+C(t^{\f12}+\v)\Lambda(\mathcal{N}(t)).
\end{align}
\end{lemma}
\noindent\textbf{Proof.}
Firstly, it follows from \eqref{2.85} and \eqref{2.94}  that 
\begin{align}\label{2.102}
\int_0^t\|(\v\partial_t)^su_x(\tau)\|^2_{L^2}d\tau\leq C\Big(1+\Lambda(\mathcal{Q}(0))\Big)+C(t^{\f12}+\v)\Lambda(\mathcal{N}(t)).
\end{align}
Using \eqref{2.102} and $\eqref{2.86}_1$, one can obtain that 
\begin{align}\label{2.103}
\int_0^t\|(\v\partial_t)^{s+1}p(\tau)\|^2_{L^2}d\tau&\leq C\int_0^t\|(\v\partial_t)^su_x(\tau)\|^2_{L^2}d\tau+ C\Big(1+\Lambda(\mathcal{Q}(0))\Big)+C(t^{\f12}+\v)\Lambda(\mathcal{N}(t))\nonumber\\
& \leq C\Big(1+\Lambda(\mathcal{Q}(0))\Big)+C(t^{\f12}+\v)\Lambda(\mathcal{N}(t)).
\end{align}

On the other hand, it follows from $\eqref{2.86}_2$ that 
\begin{align}\label{2.105}
(\v\partial_t)^sp_{x}=-\v b(-\t)(\v\partial_t)^su_t+(\v\partial_t)^s\Big\{-\v b(-\t)u u_{x}+\tilde{\mu}a(\v p)\v u_{xx}\Big\}.
\end{align}
Multiplying \eqref{2.105} by $(\v\partial_t)^sp_{x}$ and integrating the resulting equation yield that 
\begin{align}\label{2.106}
\|(\v\partial_t)^sp_{x}\|^2_{L^2}&=-\v \int b(-\t)(\v\partial_t)^su_t\cdot (\v\partial_t)^sp_{x}dx\nonumber\\
&~~~+\int (\v\partial_t)^sp_{x}\cdot (\v\partial_t)^s\Big\{-\v b(-\t)u u_{x}+\tilde{\mu}a(\v p)\v u_{xx}\Big\}dx\nonumber\\
&\leq -\f{d}{dt}\int b(-\t)(\v\partial_t)^su\cdot (\v\partial_t)^s(\v p_{x})dx+\int b(-\t)(\v\partial_t)^su\cdot (\v\partial_t)^{s+1}p_{x}dx\nonumber\\
&~~~+\int \partial_t b(-\t)(\v\partial_t)^su\cdot (\v\partial_t)^{s}(\v p)_{x}dx+\f18\|(\v\partial_t)^sp_{x}\|^2_{L^2}\nonumber\\
&~~~+\int \Big|(\v\partial_t)^s\Big\{-\v b(-\t)u u_{x}+\tilde{\mu}a(\v p)\v u_{xx}\Big\}\Big|^2dx.
\end{align}
Integrating \eqref{2.106} with respect to time,  using \eqref{2.102}, \eqref{2.103}, Lemma \ref{lem2.1}-Lemma \ref{lem2.3} and Lemma \ref{lem2.7}, one immediately gets that 
\begin{align}\label{2.107}
\int_0^t\|(\v\partial_t)^sp_{x}(\tau)\|^2_{L^2}d\tau\leq C\Big(1+\Lambda(\mathcal{Q}(0))\Big)+C(t^{\f12}+\v)\Lambda(\mathcal{N}(t)).
\end{align}
Then, from \eqref{2.86},  \eqref{2.103} and \eqref{2.107}, one obtains that 
\begin{align}\label{2.108}
&\int_0^t\|(\v\partial_t)^{s-1}(p_{xx},u_{xx})(\tau)\|^2_{L^2}d\tau
\nonumber\\
&\leq C\int_0^t\|(\v\partial_t)^{s}(p_{x},u_{x})(\tau)\|^2_{L^2}d\tau+C\Big(1+\Lambda(\mathcal{Q}(0))\Big)+C(t^{\f12}+\v)\Lambda(\mathcal{N}(t))\nonumber\\
&\leq C\Big(1+\Lambda(\mathcal{Q}(0))\Big)+C(t^{\f12}+\v)\Lambda(\mathcal{N}(t)).
\end{align}
In the same way, one can prove that 
\begin{align}
&\int_0^t\|(p_x,u_x)(\tau)\|^2_{\mathcal{H}^s}d\tau
\leq C\Big(1+\Lambda(\mathcal{Q}(0))\Big)+C(t^{\f12}+\v)\Lambda(\mathcal{N}(t)).\nonumber
\end{align}
Therefore the proof of Lemma \ref{lem2.8} is completed.
$\hfill\Box$

\

\noindent\textbf{Proof of Proposition \ref{pro2.1}}. Proposition \ref{pro2.1} follows immediately from Lemma \ref{lem2.1}-Lemma \ref{lem2.3}, Lemma \ref{lem2.7} and Lemma \ref{lem2.8}.
$\hfill\Box$

\

\noindent\textbf{Proof of Theorem \ref{thm1.1}}.  Using Proposition \ref{pro2.1}, one can prove Theorem \ref{thm1.1} by combining the local existence theorem and  the bootstrap arguments, and thus close the a priori assumption \eqref{2.14}.  The details are omitted here for simplicity of presentation.  
$\hfill\Box$

\subsection{Low Mach Limit}

In this subsection, we will prove Theorem \ref{ipC} with a modified compactness argument, which  was introduced by  M\'{e}tivier and Schochet in \cite{MS}, see also the extensions in \cite{Alazard1, Alazard2, Levermore}.

\noindent\textbf{Proof of Theorem \ref{ipC}}. Firstly, it follows from Theorem \ref{thm1.1} that
\begin{equation}\label{2.112}
\sup_{\tau\in[0, T_0]}\|(p^\v, u^\v)(\tau)\|_{H^s}+\sup_{\tau\in[0, T_0]}\|(\theta^\v-\tilde{\theta})(\tau)\|_{H^{s+1}}<\infty.
\end{equation}
Then extracting a subsequence, it holds that 
\begin{equation*}
(p^\v, u^\v)\rightharpoonup(\overline{p}, \overline{u}) ~~ \mbox{as}~~ \v\rightarrow 0 ~~\mbox{weak}-* ~~\mbox{in}~~ L^\infty(0, T_0; H^{s}(\mathbb{R})),
\end{equation*}
\begin{equation*}
\theta^\v-\tilde{\theta}\rightharpoonup \bar{\theta}-\tilde{\theta}~~ \mbox{as}~~ \v\rightarrow 0 ~~\mbox{weak}-* ~~\mbox{in}~~ L^\infty(0, T_0; H^{s+1}(\mathbb{R})).
\end{equation*}
It follows from the equation of $\theta^\v$ and \eqref{2.112} that 
$$\theta^\v_t\in L^{\infty}(0, T_0; H^{s-2}(\mathbb{R})).$$ 
which, together with Aubin-Lions lemma, yields that the functions $\theta^\v$ converge (possibly after extracting a subsequence) to $\bar{\theta}$ strongly in $C([0, T_0]; H^{s'+1}_{loc}(\mathbb{R}))$  for all $s'<s$.

To obtain the limiting system \eqref{2.8}, we need to show the strong convergence of $(p^\v, u^\v)$ in $L^2(0,T_0; H^{s'}_{loc}(\mathbb{R}))$ for  $s'<s$. To this end, we will show $p^\v$ and $(2u^\v-\kappa e^{\theta^\v-\v p^\v}\theta^\v_x)_x$ converge strongly to $0$ as $\v\rightarrow0$. In fact we rewrite $(\ref{2.1})_1$ and $(\ref{2.1})_2$ as,
\begin{equation}\label{LL1}
\v p^\v_t+(2u^\v-\kappa e^{\theta^\v-\v p^\v}\theta^\v_x)_x=\v f^\v,
\end{equation}
and 
\begin{equation}\label{LL2}
\v e^{-\theta^\v} u^\v_t+p^\v_x=\v g^\v.
\end{equation}
It follows from \eqref{ueid} that $f^\v$ and $g^\v$ are uniformly bounded in $C([0, T_0]; H^{s-1}(\mathbb{R}))$. Noticing $p^\v$ are uniformly bounded in $L^\infty(0, T_0, L^\infty(\mathbb{R}))$, passing the weak limit in \eqref{LL1} and \eqref{LL2} leads to  $(\overline{p})_x=0$ and $(2\overline{u}-\kappa e^{\bar{\theta}}\bar{\theta}_x)_x=0$ in the distribution sense.

On the other hand, by taking the limit of  $(\ref{2.1})_3$, we can get that  $\bar{\theta}$ satisfies
\begin{equation}\label{2.8-2}
\bar{\theta}_t=\f{\k e^{\bar{\theta}}}{2}\bar{\theta}_{xx},
\end{equation} 
with the initial data,
\begin{equation}\label{2.8-2i}
\bar{\t}( x, 0)=\t_{in}(x).
\end{equation} 
From the maximum principle and energy method, one can show the existence and uniqueness of smooth solution of \eqref{2.8-2}, \eqref{2.8-2i}. To get the spatial decay of $\bar{\t}$, as in \cite{Alazard2}, we define:
\begin{equation*}
H=x^{1+\sigma}(\bar{\t}-\t_+),
\end{equation*}
which satisfies
\begin{equation*}
H_t=\f{\k e^{\bar{\t}}}{2}H_{xx}-\f{\k(1+\sigma) e^{\bar{\t}}}{x}H_x+\f{\k(1+\sigma)(2+\delta) e^{\bar{\t}}}{x^2} H.
\end{equation*}
By the energy estimates and \eqref{ueid}, we can obtain
$$
\|H\|_{L^{\infty}(0, T_0; H^1[1,\infty))}\leq C[\|H(0)\|_{ H^1[1,\infty)}+\|\bar{\t}-\tilde{\t}\|_{H^2}+1]
\leq C[1+\Lambda(C_0)]
$$ 
which, together with  Sobolev embedding, yields that
\begin{equation} \label{DCC1}
|\bar{\theta}(x, t)-\theta_+|\leq C x^{-1-\sigma}, \mbox{ as } x\in [1, +\infty).
\end{equation}

To obtain the strong convergence of $u^\v$ and $p^\v$, we need the following Proposition \ref{CCL} which will be proved in the end of this section.

\begin{proposition}\label{CCL}
Let \eqref{ueid} and \eqref{DCC1} hold, then $p^\v$ and $(2u^\v-\kappa e^{\theta^\v-\v p^\v}\theta^\v_x)_x$ converge to $0$ strongly in $L^2(0, T_0; H^{s'}_{loc}(\mathbb{R}))$ and  $L^2(0, T_0; H^{s'-1}_{loc}(\mathbb{R}))$ for $s'<s$, 
respectively. 
\end{proposition}
If Proposition \ref{CCL} holds, passing the limit in the equations \eqref{2.1} for   $(p^\v, u^\v, \theta^\v)$, one proves that the limit $(0, \overline{u}, \overline{\theta})$ solves \eqref{2.8} in the sense of distribution.

On the other hand, following the arguments in \cite{MS}, one can obtain that $(\bar u,\bar{\t})$ satisfies the initial condition
\begin{equation}\label{2.121}
(\bar{u},\bar{\t})|_{t=0}=(w_{in},\t_{in})
\end{equation}
where $w_{in}$ is determined by $w_{in}=\f12\k e^{\t_{in}}(\t_{in})_x$.
Moreover one can get the uniqueness of solutions to the limit system \eqref{2.8} with initial data \eqref{2.121} by the energy method and then the above conclusions hold for the whole sequence $(p^\v, u^\v, \theta^\v)$. The proof of Theorem \ref{ipC} is completed. $\hfill\Box$

\

\noindent\textbf{Proof of Proposition \ref{CCL}:}  Applying $\v \partial_t$ to \eqref{LL1} and $\partial_x$  to $e^{\theta^\v}\times\eqref{LL2}$, we  obtain that
\begin{equation}\label{pce}
\v^2\f12 p^\v_{tt}-(e^{\theta^\v}p^\v_x)_x=\v F^\v(p^\v, u^\v, \theta^\v),
\end{equation}
where $F^\v(p^\v, u^\v, \theta^\v)$ is a smooth function. From  \eqref{ueid},  $\v F^\v(p^\v, u^\v, \theta^\v)$ converges to $0$  strongly in $L^2(0, T_0; L^2(\mathbb{R}))$ as $\v\rightarrow0$. We now recall the convergence lemma due to M\'{e}tivier and Schochet \cite{MS} in $\mathbb{R}^d$ with $d=1, 2, 3$.
\begin{lemma}\label{MSCL}
	Let $T>0$, $v^\varepsilon$ be a bounded sequence in $C([0, T], H^2(\mathbb{R}^d))$ and $\varepsilon\partial_tv^\varepsilon$ are bounded in $L^2(0, T; L^2(\mathbb{R}^d))$ satisfying
	\begin{equation}
	\varepsilon^2\partial_t(a^\varepsilon\partial_t v^\varepsilon)-\nabla \cdot (b^\varepsilon \nabla v^\varepsilon )=c^\varepsilon,
	\end{equation}
	where $c^\varepsilon$ converges to $0$ strongly in $L^2(0, T; L^2(\mathbb{R}^d))$. Assume further that, for some $k>1+d/2$, the coefficients $(a^\varepsilon, b^\varepsilon)$ are positive and uniformly bounded in $C([0, T]; H^k_{loc}(\mathbb{R}^d))$ and converge in $C([0, T]; H^k_{loc}(\mathbb{R}^d))$ to  $(a, b)$ satisfying
	\begin{equation}\label{KR}
	\mbox{for all }\tau\in\mathbb{R}, \mbox{ the kernel of } a\tau^2+\nabla \cdot (b \nabla) \mbox{ in } L^2(\mathbb{R}^d) \mbox{ is reduced to } \{0\}.
	\end{equation}
	Then the sequence $v^\varepsilon$ converges to 0 strongly in $L^2(0, T; L^2_{loc}(\mathbb{R}^d))$.
\end{lemma}

\

We introduce the following condition to verify \eqref{KR}.\\[1mm]
\noindent\textbf{Condition A:} {\it The functions $(a, b)$ are positive bounded and satisfy
	\begin{equation*}
	|a(x, t)-a_+|\leq C_{a}r(x), \quad |b(x, t)-b_+|\leq C_{b}r(x), \mbox{ as } x\in [1, +\infty).
	\end{equation*}
	where $r(x)\in L^1([1, +\infty))$ is a non-negative function, and $a_+$, $b_+$, $C_a$,  $C_b$ are some positive constants.}

\

It is obvious that if $v\in L^2(\mathbb{R})$ satisfies
\begin{equation}\label{Keq}
a\tau^2 v+\partial_x(b \partial_x v)=0,
\end{equation}
then $v\in H^1(\mathbb{R})$ and
\begin{equation*}
\int_{\mathbb{R}}b(\partial_x v)^2 dx= \tau^2 \int_{\mathbb{R}}a  v^2 dx,
\end{equation*}
which implies $v \equiv 0$ when  $\tau=0$. When $\tau\neq 0$, we assume  $\tau=1$ without loss of generality. Let $w := b \partial_x v$,  then \eqref{Keq} becomes
\begin{equation}\label{ODE1}
\frac{d}{d x}\left(\begin{array}{c}
v\\
w
\end{array} \right)
=\left(\begin{array}{cc}
0&h\\
-a&0
\end{array} \right)
\left(\begin{array}{c}
v\\
w
\end{array} \right),
\end{equation}
where $h:=1/b$ has  similar properties of $b$. From $v\in H^1(\mathbb{R})$, $v, w\rightarrow 0$ as $|x|\rightarrow \infty$. 

For further analysis, we need the following Lemma \ref{ODEC} which can be verified directly by the energy method. The details are omitted here. 
\begin{lemma}\label{ODEC}
Let $A(x)$ be a smooth bounded function on $x\in\mathbb{R}^+$. Assume $\frac{d}{d x} U= A(x) U$,  $U(0)=0$ with $U\in \mathbb{C}^d$.  Then $U \equiv 0$.
\end{lemma}

We shall first prove $(v,w)^T(x) = 0$ for  $x\in[1, +\infty)$. It is noted that  $(a, h)\rightarrow (a_+, h_+)$ as $x\rightarrow+\infty$. Then we rewrite \eqref{ODE1} as
\begin{equation}\label{2.128}
\frac{d}{d x}\left(\begin{array}{c}v\\w \end{array} \right)+\left(\begin{array}{cc}
0&-h_+\\
a_+&0
\end{array} \right)
\left(\begin{array}{c}
v\\
w
\end{array} \right)
=\left(\begin{array}{cc}
0&h-h_+\\
-(a-a_+)&0
\end{array} \right)
\left(\begin{array}{c}
v\\
w
\end{array} \right) =:\left(\begin{array}{cc}
0&\tilde{h}\\
-\tilde{a}&0
\end{array} \right)
\left(\begin{array}{c}
v\\
w
\end{array} \right).
\end{equation}
A direct calculation shows that $B=Q ^{-1}\Lambda Q$ with
\begin{equation}\label{2.129}
B=\left(\begin{array}{cc}
0&-h_+\\
a_+&0
\end{array} \right),~~
\Lambda=\left(\begin{array}{cc}
\sqrt{a_+h_+}i&0\\
0&-\sqrt{a_+h_+}i
\end{array} \right)
~~\mbox{and}~~Q=\left(\begin{array}{cc}
\sqrt{a_+}&\sqrt{h_+} i\\
\sqrt{a_+}&~-\sqrt{h_+}i
\end{array} \right).
\end{equation}
It follows from \eqref{2.128} and \eqref{2.129} that
\begin{align}\label{2.130}
&\frac{d}{d x}\left(\begin{array}{c}\sqrt{a_+} v+\sqrt{h_+} i w \\\sqrt{a_+} v-\sqrt{h_+} i w  \end{array} \right)+\left(\begin{array}{cc}
\sqrt{a_+h_+}i&0\\
0&-\sqrt{a_+h_+}i
\end{array} \right)
\left(\begin{array}{c}\sqrt{a_+} v+\sqrt{h_+} w i\\\sqrt{a_+} v-\sqrt{h_+} w i \end{array} \right)\nonumber\\
&=\left(\begin{array}{cc}
\sqrt{a_+}&\sqrt{h_+} i\\
\sqrt{a_+}&~-\sqrt{h_+}i
\end{array} \right)\left(\begin{array}{cc}
0&\tilde{h}\\
-\tilde{a}&0
\end{array} \right)
\left(\begin{array}{c}
v\\
w
\end{array} \right),
\end{align}
which yields immediately that
\begin{equation*}\label{2.131}
\frac{d}{d x}\left(\begin{array}{c}e^{\sqrt{a_+h_+}i x}\left(\sqrt{a_+} v+\sqrt{h_+} w i\right)\\e^{-\sqrt{a_+h_+}i x}\left(\sqrt{a_+} v-\sqrt{h_+} w i\right) \end{array} \right) =\left(\begin{array}{cc}
\sqrt{a_+}e^{\sqrt{a_+h_+}ix} &\sqrt{h_+} i e^{\sqrt{a_+h_+}ix} \\
\sqrt{a_+}e^{-\sqrt{a_+h_+}ix} &~-\sqrt{h_+} i e^{-\sqrt{a_+h_+}ix}
\end{array} \right)
\left(\begin{array}{cc}
0 &\tilde{h}\\
-\tilde{a} &0
\end{array} \right)
\left(\begin{array}{c}
v\\
w
\end{array} \right).
\end{equation*}
Then, one has that 
\begin{align}
	\frac{d}{d x}\left(\begin{array}{c}e^{\sqrt{a_+h_+}i x}\left(\sqrt{a_+} v+\sqrt{h_+} w i\right)\\e^{-\sqrt{a_+h_+}i x}\left(\sqrt{a_+} v-\sqrt{h_+} w i\right) \end{array} \right)=\tilde{B}
	\left(\begin{array}{c}e^{\sqrt{a_+h_+}i x}\left(\sqrt{a_+} v+\sqrt{h_+} w i\right)\\e^{-\sqrt{a_+h_+}i x}\left(\sqrt{a_+} v-\sqrt{h_+} w i\right) \end{array} \right),\nonumber
\end{align}
where
\begin{eqnarray}\label{2.133}
\tilde{B}&:=&\left(\begin{array}{cc}
\sqrt{a_+}e^{\sqrt{a_+h_+}ix}&\sqrt{h_+} i e^{\sqrt{a_+h_+}ix} \\
\sqrt{a_+}e^{-\sqrt{a_+h_+}ix}&~-\sqrt{h_+} i e^{-\sqrt{a_+h_+}ix}i
\end{array} \right)
\left(\begin{array}{cc}
0&~\tilde{h}\\
-\tilde{a}&~0
\end{array} \right)\left(\begin{array}{cc}
\sqrt{a_+}e^{\sqrt{a_+h_+}ix} &\sqrt{h_+} i e^{\sqrt{a_+h_+}ix} \\
\sqrt{a_+}e^{-\sqrt{a_+h_+}ix} &~-\sqrt{h_+} i e^{-\sqrt{a_+h_+}ix}
\end{array} \right)^{-1}\nonumber\\
&=&\frac{i}{2\sqrt{a_+h_+}}\left(\begin{array}{cc}
-h_+\tilde{a}-a_+\tilde{h} &(a_+\tilde{h}-h_+\tilde{a}) e^{2\sqrt{a_+h_+}ix} \\
(h_+\tilde{a}-a_+\tilde{h})e^{-2\sqrt{a_+h_+}ix} &h_+\tilde{a}+a_+\tilde{h}
\end{array} \right).\nonumber
\end{eqnarray}

We introduce a new coordinate $y=\int_{x}^{+\infty}r(z)dz$ satisfying
\begin{equation*}
\frac{d}{d x}= \frac{d y}{d x} \frac{d}{d y}=-r(x)  \frac{d}{d y}.
\end{equation*}
Define
\begin{equation*}\label{2.134}
A:=\frac{-1}{r(x)}\tilde{B}~~\mbox{and}~~U(y):=\left(\begin{array}{c}e^{\sqrt{a_+h_+}i x}\left(\sqrt{a_+} v+\sqrt{h_+} w i\right)\\e^{-\sqrt{a_+h_+}i x}\left(\sqrt{a_+} v-\sqrt{h_+} w i\right) \end{array} \right).
\end{equation*}
From \textbf{Condition A}, we can check that  $A$ is a smooth  bounded function. The initial data $U(0)=0$ follows from the fact $\lim_{x\rightarrow\infty}(v, w)(x)=0$.  Therefore, it follows from Lemma \ref{ODEC} that  $(v, w)(x)\equiv0$ for $x\in[1, +\infty)$. Going back to the original ODE system \eqref{ODE1} and employing Lemma \ref{ODEC} again, one can prove $v\equiv0$ from which \eqref{KR} holds.

Now, we return to the convergence of $p^\v$. By the strong convergence of $\theta^\v$ and \eqref{DCC1}, one can prove that the coefficients in \eqref{pce} satisfy  \textbf{Condition A}. It follows from  Lemma \ref{MSCL} that $p^\v\rightarrow0$ strongly in $L^2(0, T_0; L^2_{loc}(\mathbb{R}))$  as $\v\rightarrow0$. On the other hand,  the uniform boundedness \eqref{ueid} and the interpolation theorem yield immediately that
\begin{equation*}
p^\v\rightarrow 0~~\mbox{strongly in}~~L^2(0, T_0; H^{s'}_{loc}(\mathbb{R})),~\mbox{for}~s'<s.
\end{equation*}
In the same way, we can prove $(2u^\v-\kappa e^{\theta^\v-\v p^\v}\theta^\v_x)_x\rightarrow0$ as $\v\rightarrow0$. Therefore the proof of Proposition \ref{CCL} is completed. 
$\hfill\Box$

\

\noindent {\bf Acknowledgments:}
 Feimin Huang is partially supported by by National Center for Mathematics
 and Interdisciplinary Sciences, AMSS, CAS and NSFC Grant No.11371349. The research of Tian-Yi Wang was supported in part
by  National Natural Sciences Foundation of China No. 11601401. Yong Wang is partially supported by National Natural
Sciences Foundation of China No. 11401565.


\begin{thebibliography}{99}

	\bibitem{Alazard1} 
	\newblock T. Alazard,  
	\newblock \emph{ Incompressible limit of the nonisentropic Euler equations with the solid wall boundary conditions} 
	\newblock Advances in Differential Equations, \textbf{10(1)} (2005), 19-44.

	
	\bibitem{Alazard2}
	\newblock T. Alazard,
	\newblock\emph{ Low Mach number limit of the full Navier-Stokes equations}. 
	\newblock Arch. Rat. Mech. Anal., \textbf{180(1)} (2006), 1-73.
	
	\bibitem{Asano}
	\newblock K. Asano,
	\newblock \emph{On the incompressible limit of the compressible Euler equations,}
	\newblock Japan J. Appl. Math., \textbf{4} (1987) 455-488.

	\bibitem{Atkinson-Peletier}
	\newblock  F.V. Atkinson and  L.A. Peletier,
	\newblock  \emph{Similarity solutions of the nonlinear diffusion equation},
	\newblock  Arch. Rat. Mech. Anal., \textbf{54} (1974), 373-392.	
	
	
	\bibitem{Batchelor} G. K. Batchelor,
	\newblock ``An introduction to fluid dynamics",
	\newblock Cambridge University Press, 2000.
	
	\bibitem{BDG} D. Bresch, B. Desjardins and E. Grenier, 
	\newblock \emph{Oscillatory limit with changing eigen values: a formal study}, in: A. V. Fursikov, G. P. Galdi and V. V. Pukhnachev (Eds.),``New Directions in Mathematical Fluid Mechanics", 
	\newblock Birkh\"{a}user, Basel, 2010, pp.91-105.
	
	\bibitem{DBDL}
	\newblock D. Bresch, B. Desjardins, E. Grenier and C. -K. Lin, 
	\newblock \emph {Low Mach number limit of viscous polytropic flows: formal asymptotics in the periodic case,}
	\newblock Stud. Appl. Math., \textbf{109} (2002) 125-149.
	
	
	\bibitem{Danchin}
	\newblock R. Danchin,
	 \newblock \emph{Low Mach number limit for viscous compressible flows,} 
	 \newblock ESAIM: Mathematical Modelling and Numerical Analysis, \textbf{39} (2005) 459-475.

	\bibitem{DeG}
	\newblock B. Desjardins and E. Grenier,
	\newblock \emph{Low Mach number limit of viscous compressible flows in the whole space,}  
	\newblock Proc. R. Soc. Lond. Ser. A Math. Phys. Eng. Sci., \textbf{455} (1999) 2271-2279.
	
	\bibitem{DGLM}
	\newblock B. Desjardins, E. Grenier, P.-L. Lions and N. Masmoudi,
	\newblock \emph {Incompressible limit for solutions of the isentropic Navier-Stokes equations with Dirichlet boundary conditions,} 
	\newblock J. Math. Pures Appl., \textbf{78} (1999) 461-471.
	
	\bibitem{Dou-Jiang-Ou}
	\newblock C. Dou, S. Jiang and Y. Ou
	\newblock \emph{Low Mach number limit of full Navier-Stokes equations in a 3D bounded domain}
	\newblock J. Differential Equations, \textbf{258} (2015) 379-398.
	
		
		\bibitem{Duyn-Peletier}
		\newblock C.T. Duyn and L.A. Peletier,
		\newblock \emph{A class of similarity
			solution of the nonlinear diffusion equation},
		\newblock Nonlinear Analysis, T.M.A., \textbf{1} (1977), 223-233.
	
	

	
	
	
	
	
	

	

	
	\bibitem{Fan-Gao-Guo}
	\newblock J. Fan, H. Gao and B. Guo,
	\newblock \emph{Low Mach number limit of the compressible magnetohydrodynamic equations with zero thermal conductivity coefficient,}
	\newblock Math. Methods Appl. Sci. \textbf{34} (2011) 2181-2188.
	
	\bibitem{Feireisl-N} 
	\newblock E. Feireisl, A. Novotny,
	\newblock `` Singular Limits in Thermodynamics of Viscous Fluids",
	\newblock Birkh\"{a}user, Basel, 2009.
	

	
	

	
	\bibitem{Hu-Wang}
	\newblock X.P. Hu and D.H. Wang,
	\newblock \emph{Low Mach number limit of viscous compressible magnetohydrodynamic flows,}
	\newblock SIAM J. Math. Anal., \textbf{41} (2009) 1272-1294.
	
	
	
	
	\bibitem{Huang}
	\newblock F. M. Huang,  
	\newblock \emph{Thermal creep flow for the Boltzmann equation.},
	\newblock Chin. Ann. Math. Ser. B \textbf{36} (2015), no. 5, 855-870. 
	
	
	\bibitem{Huang-Wang-Wang-Yang-2} 
	\newblock F. M. Huang, Y. Wang, Y. Wang and T. Yang,
	\newblock \emph{Justification of Diffusion limit for the Boltzmann Equation with a non-trivial Profile},
	\newblock  Quart. Appl. Math. \textbf{74} (2016), no. 4, 719-764. 
	
	\bibitem{Huang-2005} 
	\newblock F.M. Huang, A. Matsumura, Z.P. Xin, 
	\newblock \emph{Stability of contact discontinuities for the 1-D compressible Navier-Stokes equations.}
	\newblock Arch. Ration. Mech. Anal., \textbf{179} (2006), 55-77.
	
	\bibitem{Iguchi}
	\newblock T. Iguchi,
	\newblock \emph{The incompressible limit and the initial layer of the compressible Euler equation in $\mathbb{R}_+^n$,}
	\newblock Math. Methods Appl. Sci., \textbf{20} (1997) 945-958.
	
	\bibitem{Isozaki1}
	\newblock H. Isozaki,
	\newblock \emph{Singular limits for the compressible Euler equation in an exterior domain,}
	\newblock J. Reine Angew. Math., \textbf{381} (1987) 1-36.
	
	\bibitem{Isozaki2}
	\newblock H. Isozaki, 
	\newblock \emph{Singular limits for the compressible Euler equation in an exterior domain, I. Bodies in an uniform flow,} \newblock Osaka J. Math., \textbf{26} (1989) 399-410.
	
	\bibitem{JJL1}
	\newblock S. Jiang, Q.C. Ju and F.C. Li,
	\newblock \emph{Incompressible limit of the compressible magnetohydrodynamic equations with vanishing viscosity coefficients,}
	\newblock SIAM J. Math. Anal., \textbf{42} (2010) 2539-2553.
	
	
	\bibitem{JJL2}
	\newblock S. Jiang, Q.C. Ju and F.C. Li,
	\newblock \emph{Low Mach number limit for the multi-dimensional full magnetohydrodynamic equations,}
	\newblock Nonlinearity \textbf{15} (2012) 1351-1365.
	
	
	\bibitem{JJL3}
	\newblock S. Jiang, Q.C. Ju and F.C. Li,
	\newblock \emph{Incompressible limit of the non-isentropic ideal magnetohydrodynamic equations,}
	\newblock SIAM J. Math. Anal. \textbf{48} (2016), no. 1, 302-319. 
	
		
		\bibitem{JJLX}
		\newblock S. Jiang, Q.C. Ju, F. Li, and Z.P. Xin,
		\newblock \emph{Low Mach number limit for the full compressible magnetohydrodynamic equations with general initial data,} 
		\newblock Advances in Mathematics, \textbf{259} (2014), 384-420.
	
	\bibitem{Jiang-Ou} 
	\newblock S. Jiang and Y.B. Ou, \newblock \emph{Incompressible limit of the non-isentropic Navier–Stokes equations with well-prepared initial data in three-dimensional bounded domains,}
	\newblock J. Math. Pures Appl., \textbf{96} (2011) 1-28.
	
	\bibitem{Kim}
	\newblock H. Kim and J. Lee, 
	\emph{The incompressible limits of viscous polytropic fluids with zero thermal conductivity coefficient,}
	\newblock Comm. Partial Differential Equations, \textbf{30} (2005) 1169-1189.
	
\bibitem{Slemrod2} Y.J. Kim, M.G. Lee, M. Slemrod, Marshall,  
\emph{Thermal creep of a rarefied gas on the basis of non-linear Korteweg-theory.}
\newblock Arch. Ration. Mech. Anal. 
\textbf{215} (2015), no. 2, 353-379.

	\bibitem{KM1} 
	\newblock S. Klainerman and A. Majda, 
	\newblock \emph{Singular limits of quasilinear hyperbolic systems with
	large parameters and the incompressible limit of compressible fluids}
    \newblock Comm. Pure Appl. Math., \textbf{34} (1981) 481-524.
	
	\bibitem{KM2} 
	\newblock S. Klainerman and A. Majda, 
	\newblock \emph{Compressible and incompressible fluids}
	\newblock Comm. Pure Appl. Math.,\textbf{35} (1982), 629-653.
	
	
	
	\bibitem{Levermore}
	\newblock C. D. Levermore, W. Sun and K. Trivisa,
	\newblock \emph{A low Mach number limit of a dispersive Navier-Stokes system,}
	\newblock SIAM J. Math. Anal., \textbf{44} (2012) 176-1807.
	
	\bibitem{Lions-Masmoudi} 
	\newblock P.-L. Lions and N. Masmoudi,
	\newblock \emph{Incompressible limit for a viscous compressible fluid,}
	\newblock J. Math. Pures Appl., \textbf{77} (1998) 585-627.
	
	\bibitem{Masmoudi}
	\newblock N. Masmoudi,
	\newblock \emph {Incompressible, inviscid limit of the compressible Navier-Stokes system,} 
	\newblock Ann. Inst. H. Poincar\'{e} Anal. Non Lin\'{e}aire, \textbf{18} (2) (2001) 199-224.
	
	\bibitem{Masmoudi1}
	\newblock N. Masmoudi, 
	\newblock \emph {Examples of singular limits in hydrodynamics in :Handbook of Differential Equations: Evolutionary Equations,}
    \newblock	vol.III, Elsevier/North-Holland, Amsterdam, 2007, pp.195-275.
	
	
    \bibitem{MS}
    \newblock G. M\'{e}tivier and S. Schochet,
    \newblock \emph{The incompressible limit of the non-isentropic Euler equations,}
    \newblock Arch. Ration. Mech. Anal. \textbf{158} (2001) 61-90.
	
	\bibitem{MS1}
	\newblock G. M\'{e}tivier and S. Schochet, 
	\newblock \emph{Averaging theorems for conservative systems and the weakly compressible Euler equations,}
	\newblock J. Differential Equations, \textbf{187} (2003) 106-183.
	
	\bibitem{Schochet1}
	\newblock S. Schochet,
	\newblock \emph{The compressible Euler equations in a bounded domain: existence of solutions and the incompressible limit,}
	\newblock Comm. Math. Phys., \textbf{10} 4 (1986) 49-75.
	
	\bibitem{Schochet}
	\newblock S. Schochet,
	\newblock \emph{Fast singular limits of hyperbolic PDEs,}
	\newblock J. Differential Equations \textbf{114} (1994) 476-512.
	
		

	\bibitem{Schochet2}
	\newblock S. Schochet, 
	\newblock \emph{The mathematical theory of the incompressible limit in fluid dynamics,} in: Handbook of Mathematical Fluid Dynamics, vol. IV, 
	\newblock Elsevier/North-Holland, Amsterdam, 2007, pp. 123-157.
	
	\bibitem{Serrin}
	J. Serrin,
	\emph{Mathematical principles of classical fluid mechanics,} in: Fluid Dynamics I/Str\"{o}mungsmechanik 
	\newblock Springer Berlin Heidelberg, 1959, pp. 125-263.
	
	\bibitem{Slemrod}
	\newblock M. Slemrod,
	\newblock \emph{The Problem with Hilbert’s 6th Problem}
	\newblock Math. Model. Nat. Phenom., \textbf{10} (2015) 6-15.
	
	
	\bibitem{Ukai} 
	\newblock S. Ukai,
	\newblock \emph{The incompressible limit and the initial layer of the compressible Euler equation,}
	\newblock J. Math. Kyoto Univ., \textbf{26} (1986) 323-331.
	
\end{thebibliography}
\end{document}